\documentclass[reqno,11pt]{amsart}
\usepackage{amssymb}
\usepackage{verbatim}
\newif\ifpdf
\ifpdf
  \usepackage[pdftex]{graphicx}
  \usepackage[pdftex]{hyperref}
\else
  \usepackage{graphicx}
\fi
\numberwithin{equation}{section}       % Number formulas within sections
\theoremstyle{plain}
\newtheorem{Thm}{Theorem}[section]
\newtheorem{Prop}[Thm]{Proposition}
\newtheorem{Lemma}[Thm]{Lemma}
\newtheorem{Cor}[Thm]{Corollary}
\newtheorem{Prop-def}[Thm]{Proposition-Definition}

\theoremstyle{definition}
\newtheorem{Example}[Thm]{Example}
\newtheorem{Def}[Thm]{Definition}
\newtheorem{Remark}[Thm]{Remark}

\newcommand{\C}{{\mathbf{C}}}

\newcommand{\N}{{\mathbf{N}}}
\newcommand{\Q}{{\mathbf{Q}}}
\newcommand{\R}{{\mathbf{R}}}

\newcommand{\Rbar}{\overline{\R}_+}

\newcommand{\fm}{{\mathfrak{m}}}

\newcommand{\cA}{{\mathcal{A}}}

\newcommand{\cC}{{\mathcal{C}}}
\newcommand{\cD}{{\mathcal{D}}}

\newcommand{\cF}{{\mathcal{F}}}
\newcommand{\cH}{{\mathcal{H}}}

\newcommand{\cO}{{\mathcal{O}}}
\newcommand{\cS}{{\mathcal{S}}}
\newcommand{\cT}{{\mathcal{T}}}
\newcommand{\cV}{{\mathcal{V}}}
\newcommand{\cM}{{\mathcal{M}}}

\newcommand{\cP}{{\mathcal{P}}}

\newcommand{\hb}{{\hat{\beta}}}

\newcommand{\hphi}{{\hat{\phi}}}

\renewcommand{\th}{{\tilde{h}}}

\newcommand{\tu}{{\tilde{u}}}

\renewcommand{\=}{ : = }
\renewcommand{\a}{\alpha}

\newcommand{\e}{\varepsilon}

\newcommand{\Lloc}{L^1_{\mathrm{loc}}}
\newcommand{\Holder}{H{\"o}lder\ }
\newcommand{\eg}{e.g.\ }
\newcommand{\ie}{i.e.\ }

\newcommand{\vv}{{\vec{v}}}
\newcommand{\qand}{{\quad\text{and}\quad}}

\newcommand{\cVqm}{{\mathcal{V}_{\mathrm{qm}}}}

\newcommand{\Gast}{\Gamma^\ast}

\newcommand{\nuL}{\nu^\mathrm{L}}

\renewcommand{\div}{\operatorname{div}}

\newcommand{\vol}{\operatorname{Vol}}

\newcommand{\supp}{\operatorname{supp}}
\newcommand{\mass}{\operatorname{mass}}

\begin{document}
%
%%%%%%%%%%%%%%%%%%%%%%%%%%%%%%%%%%%%%%%%%%%%%%%%%%%%%%%%%%%%%%%%%%
%

\setcounter{tocdepth}{1}

\title{Valuative analysis of planar pluri\-sub\-harmonic functions}
\date{\today}
\author{Charles Favre \and Mattias Jonsson}
\address{CNRS-Universit{\'e} Paris 7\\
         Institut de Math{\'e}matiques\\
         Equipe G{\'e}om{\'e}trie et Dynamique\\
         F-75251 Paris Cedex 05\\
         France}
\email{favre@math.jussieu.fr}
\address{Department of Mathematics\\ 
         University of Michigan\\
         Ann Arbor, MI 48109-1109\\
         USA}
\email{mattiasj@umich.edu}
\address{Department of Mathematics\\ 
         Royal Institute of Technology\\
         SE-100 44 Stockholm\\
         Sweden}
\email{mattiasj@kth.se}
\thanks{Second author supported by NSF Grant No DMS-0200614.}


%
%%%%%%%%%%%%%%%%%%%%%%%%%%%%%%%%%%%%%%%%%%%%%%%%%%%%%%%%%%%%%%%%%%
%
\begin{abstract}
  We show that valuations on the ring $R$ of holomorphic germs
  in dimension 2 may be naturally evaluated on plurisubharmonic
  functions, giving rise to generalized Lelong numbers in the sense of
  Demailly.  Any plurisubharmonic function thus defines a real-valued
  function on the set $\cV$ of valuations on $R$ and---by way of a
  natural Laplace operator defined in terms of the tree structure on
  $\cV$---a positive measure on $\cV$.  This measure contains a great
  deal of information on the singularity at the origin. 
  Under mild regularity assumptions, it yields an exact formula for
  the mixed Monge-Amp{\`e}re mass of two 
  plurisubharmonic functions.
  As a consequence, any generalized Lelong number 
  can be interpreted as an average of valuations.  
  Using our machinery we also show that the singularity of any 
  positive closed $(1,1)$ current $T$ can
  be attenuated in the following sense: there exists a finite
  composition of blowups such that the pull-back of $T$ decomposes
  into two parts, the first associated to a divisor 
  with normal crossing support, the second having small Lelong numbers.
\end{abstract}
%
%%%%%%%%%%%%%%%%%%%%%%%%%%%%%%%%%%%%%%%%%%%%%%%%%%%%%%%%%%%%%%%%%%
%
\maketitle
%
%%%%%%%%%%%%%%%%%%%%%%%%%%%%%%%%%%%%%%%%%%%%%%%%%%%%%%%%%%%%%%%%%%
%
\tableofcontents
%
%
%%%%%%%%%%%%%%%%%%%%%%%%%%%%%%%%%%%%%%%%%%%%%%%%%%%%%%%%%%%%%%%%%%
%
%
\newpage
%
%
%%%%%%%%%%%%%%%%%%%%%%%%%%%%%%%%%%%%%%%%%%%%%%%%%%%%%%%%%%%%%%%%%%
%
%
\section*{Introduction}
Valuation theory was a fundamental tool in the work of Zariski on
desingularization of algebraic varieties. The key role in Zariski's
approach was played by what is now called the 
\emph{Riemann-Zariski surface} introduced in~\cite{Zar}. 
It is by definition a collection of valuations. 
Our aim is to use a version of the Riemann-Zariski surface 
of the ring $R$ of holomorphic germs in dimension two to
study local singularities of plurisubharmonic (psh) functions.
Our thesis, as put forward in this paper and a forthcoming 
one~\cite{FJ-mult}, is that valuations capture 
essentially all the information on a local singularity 
of a psh function in the plane. 

\smallskip
Our study is local. We say that a psh function $u$,
defined near the origin in $\C^2$, has a \emph{singularity} 
(at the origin) when $u(0)=-\infty$. To study such a singularity,
several quantities have been introduced. 
The limit of $(\log r)^{-1} \sup \{u(q) \ ; \ |q| \le r\}$ as $r\to 0$
exists~\cite{lelong} and is called the \emph{Lelong number} of $u$. 
This number gives information on the growth of $u$ around its singular
point. It coincides with the multiplicity of the curve $\{\psi=0\}$ when
$u=\log|\psi|$, $\psi\in R$. 

More generally, given 
local coordinates $(x,y)$ and weights $a,b>0$ one can~\cite{kis1}
associate to $u$ the limit of
$(ab/\log r)\sup\{u(x,y)\ ; \ |x|\le r^{1/a},\ |y|\le r^{1/b}\}$ as $r\to0$. 
This limit is known as the \emph{Kiselman number} 
and gives more precise information on the singularity.
The Kiselman number of
$u=\log|\psi|$, $\psi\in R$ coincides with $\nu(\psi)$ where $\nu$
denotes the monomial valuation sending $x$ to $a$ and $y$ to $b$.

Demailly~\cite{dem2} put these constructions into a
more general framework.  He defined a notion of psh \emph{weight} 
(see Section~\ref{sec-dem} below) and attached to any psh function $u$ and
any weight $\varphi$ a \emph{generalized Lelong number}
$\nu_\varphi(u)\ge0$. 

Three questions naturally arise:
1) Are Demailly's generalized Lelong numbers related to valuations? 
2) When do two psh weights define the same generalized Lelong number?
3) What information on the singularity of $u$ can be recovered from 
the generalized Lelong numbers?

In this paper and in~\cite{FJ-mult} we address these three questions.
As an answer to the first, we prove that any valuation may be
evaluated on psh functions (Theorem~\ref{thm-eval}), and as such
defines a generalized Lelong number (Proposition~\ref{P6}).
Conversely we prove that any generalized Lelong number associated to a
(sufficiently regular) weight is an \emph{average} of valuations
(Theorem 8.7), giving an answer to the second question. 
The averaging is with respect to a positive measure,
the \emph{tree measure} of the weight.

The tree measure $\rho_u$ of a psh function $u$ is a positive 
measure (of mass equal to the Lelong number of $u$) on a space
$\cV$ of valuations defined below. It is computed in terms of
valuations, hence is determined by the collection of all
generalized Lelong numbers of $u$. 
A slightly vague answer to the third question above is therefore 
that the tree measure $\rho_u$ 
contains essentially complete information on the singularity of $u$.
Let us outline four ways to make this assertion more precise.

First, two psh functions have the same tree measure iff 
their pullbacks by any composition of point blowups 
above the origin have 
the same Lelong numbers at any point on the exceptional divisor.

Second, we show in~\cite{FJ-mult} that the tree measure $\rho_u$
determine the multiplier ideals of all multiples of 
a psh function $u$. In that paper we also use valuations to 
give an affirmative answer to the ``openness conjecture''
by Demailly and Koll{\`a}r.

Third, we use tree measures to prove that every positive closed $(1,1)$
current $T$, defined near the origin, admits an \emph{attenuation of
  singularities}: there exists a finite composition of blowups, such
that the pullback of $T$ decomposes into two parts, the first
associated to a divisor with normal crossing support, and the second
having arbitrarily small Lelong numbers everywhere.  The corresponding
global result---for currents on compact complex surfaces---is an easy
consequence, and has recently been proved by V.~Guedj~\cite{guedj} in
a more elementary way. His proof follows the method by
Mimouni~\cite{Mim}, who was the first, to our knowledge, to study
singularities of psh functions through sequences of blowups. However,
the global result does not imply the local one, and in fact our method
gives stronger control on the second part of the decomposition, even
in the global case.

Fourth, we explore the relationship between the 
singularities of
two individual psh functions $u$ and $v$, and the mass at the origin
of the mixed Monge-Amp{\`e}re measure $dd^c u \wedge dd^c v$. For
general $u$ and $v$, only inequalities in terms of Lelong or Kiselman
numbers were known, see \eg~\cite[Chapter~3]{dem0},~\cite{rash}. Here
we obtain much sharper estimate for $dd^c u \wedge dd^c v \{0\}$ in
terms of the tree measures $\rho_u$ and $\rho_v$. 
In fact, we obtain an exact formula under a mild regularity assumption
on either $u$ or $v$. The fact that any generalized Lelong number of a
(H{\"o}lder) psh weight is an average of valuations is a consequence of
this formula. The corresponding statement for homogeneous psh
functions was obtained by Rashkovskii.

\smallskip
Let us describe more precisely the content of the article. The
Riemann-Zariski surface is classically endowed with a non-Hausdorff
topology. We described in~\cite{treeval} a natural way to 
turn it into a Hausdorff compact space, called
the \emph{valuative tree}, that is well suited for analysis. 
This space is the set $\cV$ of all normalized 
$\R_+\cup\{\infty\}$-valued valuations on $R$ 
centered at the maximal ideal. The structure
of $\cV$, essential to our analysis here,
is described in great detail in~\cite{treeval};
we recall in Section~\ref{generalities} the results
that we use in the present paper. 
Of particular importance is the fact that $\cV$ has a natural 
\emph{tree structure}: it is made up of pieces (segments) that are
canonically parameterized by real intervals.
The main difficulty in generalizing our results to higher
dimensions lies in the fact that the structure of the 
Riemann-Zariski variety is not well  understood in dimension three 
or higher.

Section~\ref{sec-psh} contains basic facts on psh functions.
In particular we recall the powerful approximation technique
due to Demailly. We also review some facts on Kiselman numbers and
(generalized) Lelong numbers.

In Section~\ref{sec-eval} we define
$\nu(u)$ for a (quasimonomial)
valuation $\nu\in\cV$ and a psh function $u$.  
The function $u\mapsto\nu(u)$ may be 
roughly characterized as the minimal, upper semicontinuous
function satisfying $\nu(\log|\psi|)=\nu(\psi)$ for all
holomorphic germs $\psi$.
Our approach is as follows.  
If $\nu$ is monomial, then we define $\nu(u)$ as a Kiselman number. 
We extend this to $\nu$ in the dense subtree $\cVqm$ of $\cV$ 
consisting of \emph{quasimonomial} valuations. 
Such valuations can be made monomial by a birational
morphism, a fact which allows us to define $\nu(u)$ 
as a growth rate of $u$ in a semi-analytic \emph{characteristic region} 
associated to $\nu$. We also identify $\nu(u)$ as a generalized 
Lelong number.

In Section~\ref{sec-further} we further investigate 
the function $u\mapsto\nu(u)$. In particular we show that
a divisorial valuation may be interpreted as the (normalized)
Lelong number at a generic point on a suitable exceptional
component.

For a fixed psh function $u$ we call the function $\nu\mapsto\nu(u)$
on $\cVqm$ the \emph{tree transform of $u$}.  It has very strong
concavity properties (some of which were noticed by Kiselman~\cite{kis1})
with respect to the tree structure on $\cV$.  
In Section~\ref{sec-pot}, we describe a class of functions
$g:\cVqm\to\R_+$ called \emph{tree potentials} and having exactly
these concavity properties.  These functions were introduced and
studied extensively in~\cite{treeval}. The set $\cP$ of tree
potentials is the smallest closed convex cone in
$\cV_\mathrm{qm}^{\R}$, containing all
functions on $\cVqm$ of the form
$\nu\mapsto\nu(\phi)$ for $\phi\in R$, 
and closed under minima. 
Theorem~\ref{riesz} asserts that there is a 
one-to-one correspondence between $\cP$ and the set
$\cM$ of positive measures on $\cV$. 
The map $\cM\to\cP$ is defined using a natural 
\emph{intersection product} on $\cV$.
Its inverse is a natural \emph{Laplace operator}
$\Delta:\cP\to\cM$.

In Section~\ref{sec-transform} we prove that the tree transform of any
psh function $u$ defines a tree potential $g_u$, hence a positive
measure $\rho_u=\Delta g_u$ on $\cV$, called the \emph{tree measure}
of $u$.  As noted above, the tree measure gives an extremely fine
description of the singularity.  Although a complete characterization
of measures arising from psh functions seems hard to obtain, we
provide some partial results.

Attenuation of singularities for currents, as described above, is
proved in Section~\ref{sec-att}.  Our proof is constructive in the
sense that the composition of blowups may be recovered from the tree
measure of the current. 
Basically, the idea is to partition the valuative
tree into subsets, each of which has small mass for the tree measure.
The subsets define points on the exceptional divisor and the mass
gives a bound for the Lelong number at the point. 
In our approach, the finiteness of a
suitable intersection product, as used in the analysis by Mimouni 
and Guedj, is replaced by the finiteness of the tree measure.

In Section~\ref{sec-inter}, we give sharp estimates of the mixed
Monge-Amp{\`e}re measure $dd^cu \wedge dd^c v \{ 0 \}$ in terms of the 
tree measures of $u$ and $v$. 
We also prove that generalized Lelong numbers
are averages of valuations, see Theorem~\ref{p-upper}. 
All results in this section rely on a reduction to the algebraic case, 
using Demailly's approximation technique.

%
%
%%%%%%%%%%%%%%%%%%%%%%%%%%%%%%%%%%%%%%%%%%%%%%%%%%%%%%%%%%%%%%%%%%
%
%
\section{The valuative tree}\label{generalities}
In this section we give a brief review of the valuative tree as
described in~\cite{treeval}. 
Throughout the paper, we set $R=\cO_0$,
the ring of holomorphic germs at the origin in $\C^2$. This is a local
ring. Its (unique) maximal ideal $\fm$ is the set of germs vanishing
at the origin, and its residue field is $\C$.
We write $(\hat{R},\hat{\fm})$ for the completion of $(R,\fm)$.
It is the ring of formal power series in two complex variables.
%
%%%%%%%%%%%%%%%%%%%%%%%%%%%%%%%%%%%%%%%%%%%%%%%%%%%%%%%%%%%%%%%%%%
%
\subsection{Valuations}
We consider the space $\cV$ of centered, normalized valuations on $R$,
\ie the set of functions $\nu:R\to[0,\infty]$ satisfying:
\begin{enumerate}
\item[(i)]\label{p1} 
  $\nu(\psi\psi')=\nu(\psi)+\nu(\psi')$ for all $\psi,\psi'$; 
\item[(ii)]\label{p2}
  $\nu(\psi+\psi')\ge\min\{\nu(\psi),\nu(\psi')\}$ for all
  $\psi,\psi'$; 
\item[(iii)]\label{p3} 
  $\nu(0)=\infty$, 
  $\nu|_{\C^*}=0$, 
  $\nu(\fm):=\min\{\nu(\psi)\ |\ \psi\in\fm\}=1$.
\end{enumerate} 
Then $\cV$ is equipped with a natural \emph{partial ordering}: 
$\nu\le\mu$ iff $\nu(\psi)\le\mu(\psi)$ for all $\psi\in\fm$. 
The \emph{multiplicity valuation} $\nu_\fm$ defined by
$\nu_\fm(\psi)=m(\psi)=\max\{k\ |\ \psi\in\fm^k\}$ 
is the unique minimal element of $\cV$. 

Note that any valuation on $R$ extends uniquely to a valuation on its
completion $\hat{R}$, hence the valuation spaces attached to
$R$ and $\hat{R}$ are isomorphic.
%
%%%%%%%%%%%%%%%%%%%%%%%%%%%%%%%%%%%%%%%%%%%%%%%%%%%%%%%%%%%%%%%%%%%%%%%%%%%
%
\subsection{Curve valuations}
Some natural maximal elements are the \emph{curve valuations} defined
as follows. To each irreducible (possibly formal) curve $C$
we associate $\nu_C\in\cV$ defined by
$\nu_C(\psi)=C\cdot\{\psi=0\}/m(C)$, where ``$\cdot$'' denotes
intersection multiplicity.
If $C$ is defined by $\phi\in\hat\fm$, then we also write 
$\nu_C=\nu_\phi$.

The set $\cC$ of local irreducible curves carries a natural
(ultra)metric in which $\cC$ has diameter 1.
It is given by $d_\cC(C,D)=m(C)m(D)/C\cdot D$.
\begin{figure}[hl]
  \includegraphics[width=0.85\textwidth]{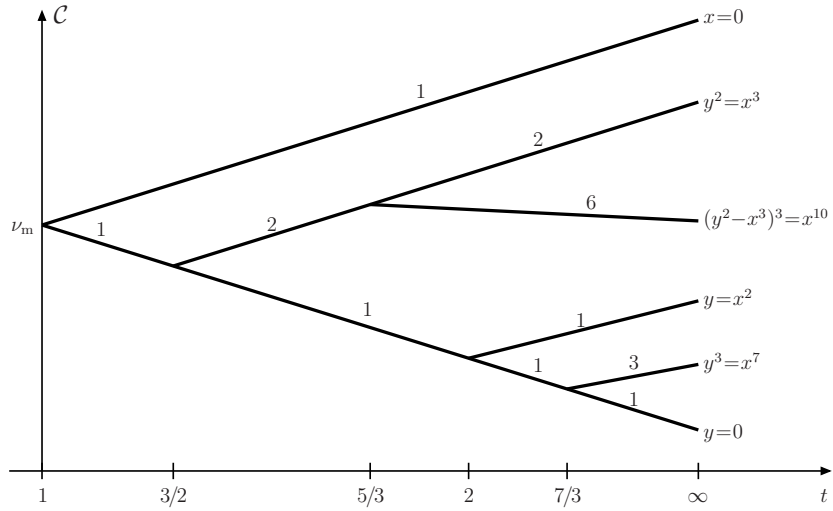}
  \caption{The valuative tree. The segments consist of valuations
    of the form $\nu_{\phi,t}$, where $\phi=x,y^2-x^3,\dots,y$
    and the skewness parameter $t$ ranges from $1$ to $\infty$.
    Skewness $t=1$ gives the multiplicity valuation $\nu_\fm$ and
    skewness $t=\infty$ the curve valuation $\nu_\phi$.
    The integer label above a segment indicates multiplicity
    (Section~\ref{sec-mult}).}\label{F4}
\end{figure}
%
%%%%%%%%%%%%%%%%%%%%%%%%%%%%%%%%%%%%%%%%%%%%%%%%%%%%%%%%%%%%%%%%%%%%%%%%%%%
%
\subsection{Quasimonomial valuations}
\begin{figure}[ht]
  \includegraphics[width=0.6\textwidth]{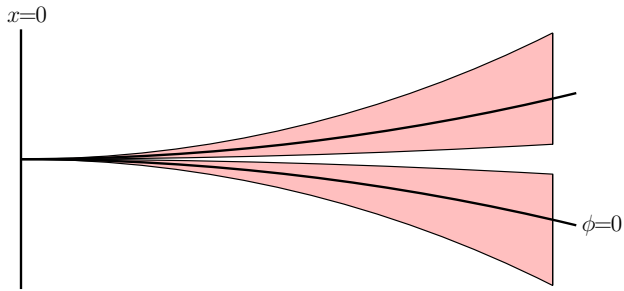}
  \caption{A characteristic region in $\C^2$.}\label{F3}
\end{figure}

Arguably the most important valuations in $\cV$ are the
\emph{quasimonomial} ones.\footnote{A quasimonomial valuation can be
  made monomial (\ie completely determined by its values on a pair of
  local coordinates $(x,y)$) by a birational morphism: see
  Section~\ref{sec-precise-mono}.  Quasimonomial valuations are also
  known as Abhyankar valuations of rank 1.}  They are of the form
$\nu_{C,t}$, where $C\in\cC$ and $t\in[1,\infty)$, and satisfy
$\nu_{C,t}(\psi)=\min\{\nu_D(\psi)\ |\ d_\cC(C,D)\le t^{-1}\}$.  We
have $\nu_{C,s}=\nu_{D,t}$ iff $s=t\ge d_\cC(C,D)^{-1}$
(see~\cite[Theorem~3.57]{treeval}).  Thus $\cVqm$, the set of all
quasimonomial valuations, is naturally a quotient of
$\cC\times[1,\infty)$, and has a natural tree structure: if
$\nu,\nu'\in\cVqm$ and $\nu<\nu'$, then the \emph{segment}
$[\nu,\nu']=\{\mu\in\cVqm\ |\ \nu\le\mu\le\nu'\}$ is isomorphic to a
compact real interval, see Figure~\ref{F4}. 
We set
$\nu_{\phi,t}:=\nu_{C,t}$ when $C=\{\phi=0\}$.  The value
$\nu_{C,t}(\psi)$ can be interpreted as the order of vanishing of
$\psi$ in a suitable ``cone'' around the curve $C$ of size
depending on $t$.
We refer to~Section~\ref{def-qm} and
Proposition~\ref{P2} for more precise statements. Figure~\ref{F3}
depicts such a cone.
Quasimonomial valuations are of two types: \emph{divisorial} and
\emph{irrational}, depending on whether the parameter $t$ is
rational or irrational.\footnote{A quasimonomial valuation $\nu$ is
  irrational iff $\nu(R)\not\subset\Q$, hence the name.}  The full
space $\cV$ is the completion of $\cVqm$ in the sense that every
element in $\cV$ is the limit of an increasing sequence in $\cVqm$.
It is hence also naturally a tree, called the \emph{valuative tree}.
The ends of $\cV$ are exactly the elements of $\cV\setminus\cVqm$
and are either curve valuations or \emph{infinitely singular}
valuations.\footnote{The latter are represented by infinite Puiseux
  series whose exponents are rational numbers with 
  unbounded denominators.  See Section~\ref{puiseux}.}
%
%%%%%%%%%%%%%%%%%%%%%%%%%%%%%%%%%%%%%%%%%%%%%%%%%%%%%%%%%%%%%%%%%%%%%%%%%%%
%
\subsection{Skewness and intersection multiplicity}\label{sec-skewness}
An important invariant of a valuation is its \emph{skewness} $\a$
defined by $\a(\nu)=\sup\{\nu(\phi)/m(\phi)\ |\ \phi\in\fm\}$.
Skewness naturally \emph{parameterizes} the trees $\cVqm$ and $\cV$ in
the sense that $\a:\cVqm\to[0,\infty)$ is strictly increasing and
restricts to a bijection onto its image on any segment; indeed
$\a(\nu_{\phi,t})=t$ for any $\nu_{\phi,t}\in\cVqm$. 
See~\cite[Section~3.3]{treeval}. 
Thus divisorial (irrational) valuations have rational
(irrational) skewness. Curve valuations have infinite skewness whereas
the skewness of an infinitely singular valuations may or may not be
finite.

The tree structure on $\cV$ implies that any collection $(\nu_i)_{i\in I}$
of valuations in $\cV$ admits an infimum $\wedge_i\nu_i$. Together with
skewness, this allows us to define an \emph{intersection product} on $\cV$:
we set $\nu\cdot\mu:=\a(\nu\wedge\mu)\in[1,\infty]$. 
This is a normalized extension of the intersection product on $\cC$ 
as $C\cdot D=(\nu_C\cdot\nu_D)m(C)m(D)$.
If $\nu\in\cV$ and $\phi\in\fm$ is irreducible, then 
$\nu(\phi)=m(\phi)\,(\nu\cdot\nu_\phi)$.
Moreover, if $\nu(\phi)$ is irrational, then 
$\nu=\nu_{\phi,t}$ with $t=\nu(\phi)/m(\phi)$.
%
%%%%%%%%%%%%%%%%%%%%%%%%%%%%%%%%%%%%%%%%%%%%%%%%%%%%%%%%%%%%%%%%%%%%%%%%%%%
%
\subsection{Tangent space and weak topology}\label{S11}
Let $\mu$ be a valuation in $\cV$. Declare 
$\nu,\nu'\in\cV\setminus\{\mu\}$ to be equivalent if
the segments $]\mu,\nu]$ and $]\mu,\nu']$ intersect.
An equivalence class is called a \emph{tangent vector} at $\mu$
and the set of tangent vectors at $\mu$, the \emph{tangent space},
denoted by $T\mu$. If $\vv$ is a tangent vector, we denote by
$U(\vv)$ the set of points in $\cV$ defining the equivalence class
$\vv$. The points in $U(\vv)$ are said to \emph{represent} $\vv$.

A point $\mu$ in the tree $\cV$ is an \emph{end}, 
a \emph{regular point},
or a \emph{branch point} when $T\mu$ contains one,
two, or three or more points, respectively.  In terms of valuations:
the ends of $\cV$ are curve and infinitely singular valuations;
the regular points irrational valuations; and the branch points
divisorial valuations, at which the tangent space is in bijection with
the complex projective line $\mathbf{P}^1$ and hence
uncountable. See~\cite[Proposition~3.20]{treeval}.

We endow $\cV$ with the \emph{weak topology},
generated by the sets $U(\vv)$ over all tangent vectors $\vv$;
this turns $\cV$ into a compact (Hausdorff) space.
The weak topology on $\cV$ is
characterized by $\nu_k\to\nu$ iff $\nu_k(\phi)\to\nu(\phi)$
for all $\phi\in R$; see~\cite[Theorem~5.1]{treeval}.
%
%%%%%%%%%%%%%%%%%%%%%%%%%%%%%%%%%%%%%%%%%%%%%%%%%%%%%%%%%%%%%%%%%%%%%%%%%%%
%
\subsection{Multiplicities}\label{sec-mult}(See~\cite[Section 3.4]{treeval})
By setting $m(\nu)=\min\{m(C)\ |\ C\in\cC,\ \nu_C\ge\nu\}$ we
extend the notion of \emph{multiplicity} from $\cC$ to $\cVqm$.  
Clearly $m:\cVqm\to\N$ is increasing and hence extends to all of $\cV$. 
In fact $m(\nu)$ divides $m(\mu)$ whenever $\nu\le\mu$.  The infinitely
singular valuations are characterized as having infinite multiplicity.

As $m$ is increasing and integer valued, it is piecewise 
constant on any segment $[\nu_\fm,\nu_\phi]$, where $\phi\in\cC$.
This implies that $m(\vv)$ is naturally defined for any tangent
vector $\vv$. If $\nu$ is nondivisorial, then $m(\vv)=m(\nu)$ 
for any $\vv\in T\nu$. 

If $\nu$ is divisorial, then the situation is more complicated.  We
refer to~\cite[Proposition~3.39]{treeval} for a precise discussion.
Suffice it to say that there exists an integer $b(\nu)$,
divisible by $m(\nu)$, such that $m(\vv)=b(\nu)$ for all
but at most two tangent vectors $\vv$ at $\nu$. We call
$b(\nu)$ the \emph{generic multiplicity} of $\nu$.
%
%%%%%%%%%%%%%%%%%%%%%%%%%%%%%%%%%%%%%%%%%%%%%%%%%%%%%%%%%%%%%%%%%%%%%%%%%%%
%
\subsection{Approximating sequences}\label{sec-approx}
(See~\cite[Section 3.5]{treeval}).
Consider a quasimonomial valuation $\nu\in\cVqm$. The
multiplicity $m$ is integer-valued and piecewise constant on the
segment $[\nu_\fm,\nu]$ and therefore has a finite number $g$
(possibly zero) of jumps.  Thus there are divisorial valuations
$\nu_i$, $0\le i\le g$ and integers $m_i$, such that
\begin{equation}\label{e705}
  \nu_\fm=\nu_0<\nu_1<\dots<\nu_g<\nu_{g+1}=\nu
\end{equation}
and $m(\mu)=m_i$ for $\mu\in\,]\nu_i,\nu_{i+1}]$, $0\le i\le g$.
We call the sequence $(\nu_i)_{i=0}^g$ the 
\emph{approximating sequence} associated to $\nu$.
It also plays a prominent role in~\cite{spiv}. 

The concept of approximating sequences extends naturally to valuations
that are not quasimonomial: for curve valuations the sequences
are still finite, for infinitely singular valuations they are infinite.
%
%%%%%%%%%%%%%%%%%%%%%%%%%%%%%%%%%%%%%%%%%%%%%%%%%%%%%%%%%%%%%%%%%%%%%%%%%%%
%
\subsection{Thinness}\label{thinness}
(See~\cite[Section 3.6]{treeval}).
Skewness $\a$ is a parameterization of $\cV$ that does not ``see'' 
multiplicities. Another parameterization, of crucial importance,
is \emph{thinness} $A$, defined as follows. 
If $\nu\in\cVqm$ then
\begin{equation}\label{e40}
  A(\nu)=2+\int_{\nu_\fm}^\nu m(\mu)\,d\alpha(\mu).
\end{equation}
In terms of~\eqref{e705} we have $A(\nu)=2+\sum_0^gm_i(\a_{i+1}-\a_i)$
with $\a_i=\a(\nu_i)$. Note that $A(\nu)\le 1+ m(\nu)\a(\nu)$.  Just
like skewness, we may define $A(\nu)$ also for $\nu\notin\cVqm$.
%
%%%%%%%%%%%%%%%%%%%%%%%%%%%%%%%%%%%%%%%%%%%%%%%%%%%%%%%%%%%%%%%%%%%%%%%%%%%
%
\subsection{Puiseux expansions}\label{puiseux}
(See~\cite[Chapter~4]{treeval}). Valuations in $\cV$ can be encoded
(nonuniquely) by Puiseux series. We describe this here in the case of
an end $\nu\in\cV$, \ie an infinitely singular or curve
valuation.  Pick local coordinates $(x,y)$ such that
$\nu_\fm<\nu\wedge\nu_y$.  Then $\nu$ is represented by
a series $\hphi=\sum_{i\ge1}a_ix^{\hb_i}$, where $a_i\in\C^*$ and
$\hb_i>1$ form an increasing sequence of rational numbers.  For
$\psi\in R$, $\psi(t,\hphi (t))=\sum b_jt^{\hat{\gamma}_j}$ is a
Puiseux series and $\nu(\psi)=\min\{\hat{\gamma}_j\ ;\ b_j\ne0\}$.
Moreover, $A(\nu)=1+\lim_{i\to\infty}\hb_i$ if the Puiseux series is
infinite and $A(\nu)=\infty$ otherwise.

If the denominators of the $\hb_i$'s are bounded, then $\nu$ is a
curve valuation. Otherwise, $\nu$ is infinitely singular.
In the latter case we can define $C_n\in\cC$ 
to be the irreducible curve associated to
the truncated Puiseux expansion $\hphi_n=\sum_1^na_ix^{\hb_i}$. 
Then $\nu_{C_n}$ converges weakly to $\nu$.
%
%%%%%%%%%%%%%%%%%%%%%%%%%%%%%%%%%%%%%%%%%%%%%%%%%%%%%%%%%%%%%%%%%%%%%%%%%%%
%
\subsection{Geometric interpretation of divisorial valuations}
\label{divisorial}
(See~\cite[Chapter~6]{treeval}). 
All divisorial valuations $\nu$ arise as follows: 
there exists a (proper) \emph{modification} 
$\pi:X\to(\C^2,0)$---in our case a finite composition of 
point blowups---and an \emph{exceptional component} $E$
(\ie $E$ is an irreducible component of the exceptional
divisor $\pi^{-1}(0)$)
such that $\nu=:\nu_E$ is equivalent to $\pi_*\div_E$, 
where $\div_E$ denotes the order of vanishing along $E$.
More precisely:
\begin{itemize}
\item
  $\nu=b^{-1}\pi_*\div_E$, where $b=b(\nu)$ is the generic
  multiplicity at $\nu$;
\item
  $A(\nu)=a/b$, where $a-1$ is equal to the order of vanishing
  along $E$ of the Jacobian determinant of $\pi$.
\end{itemize}
In fact, the pair $(a,b)$, which is called the \emph{Farey weight}
of $E$ in~\cite{treeval}, can be obtained in a purely 
combinatorial way (see also~\cite{hubbard}). 
This combinatorial procedure can be used to recover the full
tree structure (partial ordering, thinness and multiplicity)
on the valuative tree $\cV$, and the assertions above
are consequences of a much more precise result:
see~\cite[Theorem~6.22]{treeval}.

The generic multiplicity $b(\nu_E)$ is the multiplicity of any
\emph{curvette} for $\nu_E$, \ie any irreducible curve $C$ whose
strict transform under $\pi$ is smooth and intersects $E$ transversely
at a smooth point on $E$.  For any curvette we have
$\nu_C>\nu_E$. See~\cite[Section~6.6.1]{treeval}.

\begin{figure}[h]
  \includegraphics[width=0.85\textwidth]{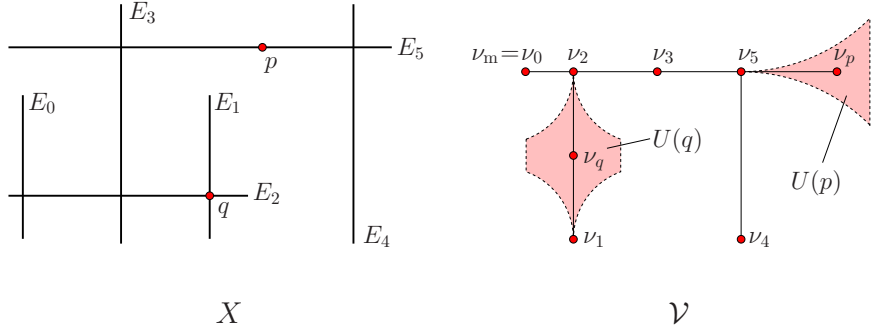}
  \caption{Dual graphs and valuations. To the left is the total space of 
a composition $\pi:X\to(\C^2,0)$ of six point blowups. 
To the right is the dual graph of $\pi$ embedded in the valuative
tree $\cV$. To the exceptional component $E_i$ is associated 
the valuation $\nu_i$. The valuations $\nu_p$ and $\nu_q$, 
obtained by blowing up the points $p$ and $q$, illustrate 
case~(i) and~(ii) of Proposition~\ref{prop-divisorial}, respectively.
The open subsets $U(p)$ and $U(q)$ of $\cV$ 
consist of all valuations in $\cV$ whose center on $X$ is
$p$ and $q$, respectively.}\label{F1}
\end{figure}

In Section~\ref{sec-att}, we shall need several results
relating divisorial valuations and their location inside $\cV$ to
dual graphs of modifications. See Figure~\ref{F1}.
\begin{Prop}{\cite[Proposition~6.37]{treeval}}\label{prop-divisorial}
  Let $\pi:X\to(\C^2,0)$ be a modification.
  Pick a point $p\in\pi^{-1}(0)$, and let $\nu_p$ 
  be the divisorial valuation associated to the blowup at $p$.
  \begin{itemize}
  \item[(i)]
    If $p$ belongs to a unique exceptional component $E$ of $\pi$, then
    $\nu_p>\nu_E$; $\nu_p$ does not represent the same tangent vector
    at $\nu_E$ as $\nu_{E'}$ for any other exceptional component 
    $E'$; and the multiplicity is 
    constant equal to $b(\nu_p)=b(\nu_E)$ on the segment $]\nu_E,\nu_p]$.
  \item[(ii)]
    If $p$ is the intersection of two exceptional components $E$ and $E'$, 
    then
    $\nu_E<\nu_p<\nu_{E'}$ (or $\nu_{E'}<\nu_p<\nu_E$); 
    and $b(\nu_p)=b(\nu_E)+b(\nu_{E'})$.
  \end{itemize}
\end{Prop}
We shall also use the following direct consequence
of~\cite[Corollary~6.39]{treeval}.
\begin{Lemma}\label{lem-estim}
  Let $\pi:X\to(\C^2,0)$ be a modification,
  and let $E,E'$ be two exceptional components 
  that intersect in $X$.
  Then $\nu_E,\nu_{E'}$ are comparable and
  \begin{equation}\label{e102}
    |\a(\nu_{E'})-\a(\nu_E)|=\frac{1}{b(\nu_{E'})b(\nu_E)}.
  \end{equation}
\end{Lemma}

%
%
%%%%%%%%%%%%%%%%%%%%%%%%%%%%%%%%%%%%%%%%%%%%%%%%%%%%%%%%%%%%%%%%%%%%%%%%%%%
%
%
\section{Plurisubharmonic functions}\label{sec-psh}
In this section we recall some facts on plurisubharmonic (psh) 
functions that we will need. General references are~\cite{dem0}
and~\cite{Hor}.
%
%%%%%%%%%%%%%%%%%%%%%%%%%%%%%%%%%%%%%%%%%%%%%%%%%%%%%%%%%%%%%%%%%%
%
\subsection{Basics}
Unless otherwise specified, all psh functions are defined in some
neighborhood of the origin in $\C^2$.
We write $dd^cu=\frac{i}{2\pi}\partial\bar\partial u$.

If $\psi$ is a holomorphic germ, \ie $\psi\in R$, then $\log|\psi|$ is psh.  
More generally, if $\psi_i\in R$, $1\le i\le n$, then
$u=\log(\sum_1^n|\psi_i|^2)$ is psh. A psh function $u$ such that 
$u-c\log(\sum_1^n|\psi_i|^2)$ is locally bounded for some 
$\psi_i\in R$ and $c>0$, is said to have \emph{logarithmic singularities}.

We shall use the following consequence of Harnack's inequality:
\begin{Lemma}\label{L-harnack}
  If $u\le0$ is psh in the bidisk $|x|<r_1$, $|y|<r_2$, then
  \begin{equation*}
    \sup_{|x|\le \e \rho_1,\ |y|\le \e \rho_2} u
    \le
    \frac{(1-\e^2)^2}{(1+\e)^4}
    \iint u(\rho_1 e^{i\theta_1},\rho_2 e^{i\theta_2})~
    \frac{d\theta_1 d\theta_2}{(2\pi)^2},
  \end{equation*}
  whenever $\rho_i\in(0,r_i)$, $i=1,2$ and $\e\in(0,1)$.
\end{Lemma}
%
%%%%%%%%%%%%%%%%%%%%%%%%%%%%%%%%%%%%%%%%%%%%%%%%%%%%%%%%%%%%%%%%%%
%
\subsection{Demailly approximation.}\label{sec-demailly}
Our study uses in a crucial 
way a powerful technique, developed by Demailly, of
approximating general psh functions by psh functions with 
logarithmic singularities.
Let us summarize the main ingredients. 
See Theorem~4.2 and its proof in~\cite{DK} for details.

Let $u$ be a psh function on a fixed ball $B$ containing the
origin. For any real number $n>0$, we let $\cH_{nu}(B)$ be the
Hilbert space of holomorphic functions $f$ on $B$ such that 
$\|f\|_{nu}:=\int_B|f|^2\exp(-2nu)\,dV<+\infty$, and define
\begin{equation*}
  u_n
  =\frac1{2n}\log\sum_1^\infty|g_{nk}|^2
  =\frac1{n}\log\sup\{|f|\ ;\ \|f\|_{nu}\le1\},
\end{equation*}
where $g_{nk}$ is an orthonormal basis of $\cH_{nu}(B)$.
Then $u_n$ converges to $u$ in a rather strong sense:
there exists $C>0$ such that for any $r>0$ small enough
\begin{equation}\label{e801}
  u(p)-\frac{C}{n}
  \le u_n(p) 
  \le \sup_{B(p,r)}u+\frac{C}{n}-\frac2n\log r.
\end{equation}
In particular this implies $u_n\to u$ in $\Lloc$.
It is the left hand inequality of~\eqref{e801} that is the harder 
one to establish---it is a consequence of the Ohsawa-Takegoshi 
extension theorem.

It is also convenient to truncate the infinite sum defining $u_n$.
More precisely, given $n$ and a smaller ball
$B'\Subset B$ containing the origin, there exists 
$k_0=k_0(u,n,B')<\infty$ such that
$\sum_1^\infty|g_{nk}|^2\le C\sum_1^{k_0}|g_{nk}|^2$
on $B'$, for some $C>0$. 
Thus $u_n-(2n)^{-1}\log\sum_1^{k_0}|g_{nk}|^2$ is bounded, 
so $u_n$ has logarithmic singularities. 

In Section~\ref{sec-inter} we will work with the class
of psh functions $u$ for which $e^u$ is \Holder continuous. 
This class, which appears at several points in the work of Demailly 
(see \eg~\cite[Theorem~2.5]{DK}), 
contains all psh functions with logarithmic singularities, but not,
for instance, the function $u(p)=-\sqrt{-\log\|p\|}$.

One consequence of \Holder continuity is that the right hand side
of~\eqref{e801} can be improved. 
Assume $u(0)=-\infty$ and that $e^u$ is \Holder continuous.
Then we may find $c>0$ so that 
$|e^{u(p)}-e^{u(p')}|\le\|p-p'\|^c$ for $p,p'$ in a
small ball around the origin. 
Thus $\sup_{B(p,r)}u\le\log(e^{u(p)}+r^c)$
for any $p$ and $r>0$ small enough. 
When $u(p)=-\infty$,~\eqref{e801} yields $u_n(p)=-\infty$ 
for $n>2/c$ by letting $r\to 0$.
If $u(p)>-\infty$, then choosing $r=e^{u(p)/c}$ 
we infer from~\eqref{e801} the estimate
\begin{equation}\label{e802}
  u-\frac{C}{n}
  \le u_n
  \le\left(1-\frac{2}{nc}\right)u+\frac{C}{n},
\end{equation}
on some neighborhood $B'$ of the origin and some constant $C>0$
(both independent of $n$).

%
%%%%%%%%%%%%%%%%%%%%%%%%%%%%%%%%%%%%%%%%%%%%%%%%%%%%%%%%%%%%%%%%%%
%
\subsection{Kiselman numbers}\label{sec-kiselman}
As we will see, the evaluation of a quasimonomial valuation on psh 
functions can be viewed as the pushforward under a birational morphism 
of a Kiselman number, so we recall some properties of the latter.
The original references are~\cite{kis1,kis2}; 
see also~\cite[Chapter~III, Example~6.11]{dem0}.

Consider a psh function $u$, fix local coordinates $(x,y)$ and
real numbers $a,b>0$. 
The \emph{Kiselman number} of $u$ with weights $(a,b)$ is defined as
\begin{equation}\label{e905}
  \nu^{x,y}_{a,b}(u)
  :=\lim_{r\to0}\tfrac{ab}{\log r}
  \sup\{u\ ;\ |x|<r^{\frac1a}, |y|<r^{\frac1b}\}.
\end{equation}
This limit exists. By the maximum principle we could equivalently
take the supremum over either the torus $\{|x|=r^{1/a},|y|=r^{1/b}\}$ 
or over the open set $\{|x|<r^{1/a}, |y|<|x|^{a/b}\}$. 
Another definition of the Kiselman number, the equivalence of which to
the previous definitions follows from Lemma~\ref{L-harnack},
is given by
\begin{equation}\label{e908}
  \nu^{x,y}_{a,b}(u)
  =\lim_{r\to0}\frac{ab}{\log r} 
  \iint u(r^{\frac1a}e^{i\theta},r^{\frac1b}e^{i\vartheta})\,
  \frac{d\theta d\vartheta}{(2\pi)^2}.
\end{equation}
In both regions $\{|x|<r^{1/a}, |y|<r^{1/b}\}$ and
$\{|x|<r^{1/a}, |y|<|x|^{a/b}\}$ we have
\begin{equation}\label{e902}
  u(p)\le\tfrac1{\min\{a,b\}}\,\nu^{x,y}_{a,b}(u)\log\|p\|+O(1)
  \quad\text{as $p\to0$}.
\end{equation}

In the particular case $a=b=1$, the Kiselman number reduces to
the \emph{Lelong number} of $u$, originally defined in~\cite{lelong}.
We denote it by $\nuL(u)$. 
The Lelong number of the positive closed current $T=dd^cu$
can be also computed as follows. 
Define $\omega=dd^c(|x|^2+|y|^2)$, and let
$\|T\|\=T\wedge\omega$ be the \emph{trace measure} of $T$. 
Then
\begin{equation}\label{eq-def-lelong}
  \nuL(T)=\lim_{r\to0}\tfrac{1}{\pi r^2}\,\|T\|\left[B(0,r)\right], 
\end{equation}
where $B(0,r)$ denotes the ball of center $0$ and radius $r$.
%
%%%%%%%%%%%%%%%%%%%%%%%%%%%%%%%%%%%%%%%%%%%%%%%%%%%%%%%%%%%%%%%%%%
%
\subsection{Intersection of currents}\label{sec-interdef}
Let us recall a few facts about intersection of currents.
See~\cite[Section~4.1]{dujardin} for more details.  Pick two psh
functions $u,v$.  When $u$ is locally integrable with respect to the
trace measure of $dd^cv$, we say that the wedge product $dd^cu\wedge
dd^cv$ is \emph{admissible}, and we define $dd^cu\wedge
dd^cv=dd^c(udd^cv)$.  In dimension two, it is a fact proved by Sibony,
see~\cite{dujardin}, that if $(u_n)_1^\infty$ and $(v_n)_1^\infty$ are
sequences of psh functions such that $u_n\ge u$, $v_n\ge v$, $u_n\to
u$ and $v_n\to v$ in $\Lloc$, then $dd^cu_n\wedge dd^cv_n$ is
admissible for all $n$ and
\begin{equation}\label{e911}
  dd^cu_n\wedge dd^cv_n\longrightarrow dd^cu\wedge dd^cv
\end{equation}
weakly, as $n\to\infty$.
As follows from~\eqref{e801}, we may apply~\eqref{e911} 
when $u_n$ and $v_n$ are Demailly approximants of
$u$ and $v$, respectively.

Equation~\eqref{e911} may also be used to show that the condition of 
$dd^cu\wedge dd^cv$ being admissible is symmetric in $u$ and $v$.
%
%%%%%%%%%%%%%%%%%%%%%%%%%%%%%%%%%%%%%%%%%%%%%%%%%%%%%%%%%%%%%%%%%%
%
\subsection{Generalized Lelong numbers}\label{sec-dem}
We now introduce some terminology and definitions taken
from~\cite[Chapter III]{dem0}
(see also the original article~\cite{dem2}). 
A psh function $\varphi$ such that
$\varphi^{-1}\{ -\infty\}=\{0\}$, and $e^{\varphi}$ is continuous is
called a \emph{psh weight}. 
If $\varphi$ is a psh weight, then $dd^c\varphi$ cannot 
charge any analytic curve. As a partial converse, 
any psh function with logarithmic singularities
and which does not charge any analytic curves defines a psh 
weight.

For any psh weight $\varphi$ and any psh function $u$, the wedge
product $dd^cu\wedge dd^c\varphi$ is admissible~\cite{dem0}, and
defines a positive measure.  We may hence associate to any psh
function $u$ a
\emph{generalized Lelong number} by setting
\begin{equation}\label{def-DL}
  \nu_\varphi (u) 
  =\lim_{r\to 0}\int_{\varphi\le\log r}dd^cu\wedge dd^c\varphi 
  =dd^cu\wedge dd^c\varphi\,\{0\}.
\end{equation}
These  generalized Lelong numbers were designed to give information 
on the nature of the singularity of $u$ at the origin. 
We will use two results by Demailly.
\begin{Prop}\cite[Chapter~III, Proposition~5.12]{dem0}.\label{p-semi}
  For any psh weight $\varphi$, and any sequence of psh functions
  $u_n\to u$ in $\Lloc$ we have 
  $\limsup\nu_\varphi(u_n)\le\nu_\varphi(u)$.
\end{Prop}
\begin{Prop}\cite[Chapter~III, Theorem~7.1]{dem0}.\label{p-comparison}
  Let $\varphi_1,\varphi_2$ be two psh weights such that
  $\limsup_{p\to0}\varphi_1(p)/\varphi_2(p)\le1$. 
  Then $\nu_{\varphi_1}(u)\le\nu_{\varphi_2}(u)$ 
  for any psh function $u$.
  In particular, equality holds when 
  $\lim_{p\to0}\varphi_1(p)/\varphi_2(p)=1$.
\end{Prop}
When the Monge-Amp{\`e}re measure $(dd^c)^2 \varphi$ is concentrated at
the origin, a Jensen type formula holds. Namely, the generalized
Lelong number of $u$ may be computed as follows 
\cite[Chapter~III, Proposition~6.5]{dem0}:
\begin{equation}\label{DL-jensen}
  \nu_\varphi(u)
  =\lim_{r\to 0}\frac{1}{\log r}\int_{\varphi=\log r}u\,d\lambda_r,
\end{equation}
where the measure $\lambda_r=(dd^c)^2\max\{\varphi,\log r\}$ 
is supported on the set $\{\varphi=\log r\}$.  

An important class of weights for which~\eqref{DL-jensen} applies is
$\varphi=\log\max\{|\phi|^s,|\psi|^t\}$ for holomorphic 
germs $\phi,\psi$ without common factor, and $s,t>0$.
The measure $\lambda_r$ is then supported on the set
$\{|\phi|^s=|\psi|^t=r\}$.

When $\varphi=\log\max\{|x|^a,|y|^b\}$ we recover the
Kiselman number:
\begin{align}
  \nu^{x,y}_{a,b}(u)
  &=dd^cu\wedge dd^c\log\max\{|x|^{a},|y|^{b}\}\{0\}\label{e906}\\
  &=\lim_{r\to0}\tfrac{1}{\log r}
  \int u\,(dd^c)^2\log\max\{r,|x|^{a},|y|^{b}\}.\notag
\end{align}
This follows from~\eqref{DL-jensen},~\eqref{e908}
and Lemma~\ref{L300} below.

We conclude this section with the following useful result.
\begin{Lemma}\label{L300}
  Fix local coordinates $(x,y)$ and real numbers $a,b,c,d\ge0$
  such that $ad\ne bc$. Let $(s,t)$ be the solution to
  $as+bt=cs+dt=1$ and assume that $s,t>0$.
  Then the measure $\lambda_r=(dd^c)^2\log\max\{r,|x|^a|y|^b,|x|^c|y|^d\}$ 
  has mass $|ad-bc|$ and 
  is proportional to the Haar-Lebesgue measure on the torus 
  $|x|=r^s$, $|y|=r^t$.
\end{Lemma}
\begin{proof}
  The measure $\lambda_r$ is supported where $|x|^a|y|^b=|x|^c|y|^d=r$
  and is invariant under rotations in $x$ and $y$, hence must be
  proportional to the Haar-Lebesgue measure on the above-mentioned
  torus. Let us prove that the mass is $|ad-bc|$.  In the case
  $a=d=1$, $b=c=0$, it is a simple application of Stokes' formula that
  the mass equals $1$. When $a,b,c,d$ are integers, consider the
  monomial change of coordinates $\pi(x,y)= (x^ay^b,x^cy^d)$. The
  topological degree of $\pi$ is equal to $|ad-bc|$. By the change of
  variable formula the mass of $(dd^c)^2\log\max\{r,|x|,|y|\}\circ \pi
  = (dd^c)^2\log\max\{r,|x|^a|y|^b,|x|^c|y|^d\}$ is $|ad-bc|$ times
  the mass of $(dd^c)^2\log\max\{r,|x|,|y|\}$, hence the mass of
  $\lambda_r$ equals $|ad-bc|$. When $a,b,c,d$ are rational numbers
  with a common denominator $q$, we use the change of variables
  $\pi(x,y) = (x^q,y^q)$ to reduce the proof to the preceding case.
  Finally we conclude in the general case by continuity.
\end{proof}
%
%
%%%%%%%%%%%%%%%%%%%%%%%%%%%%%%%%%%%%%%%%%%%%%%%%%%%%%%%%%%%%%%%%%%%%%%%%%%%
%
%
\section{Evaluating quasimonomial valuations on psh functions}\label{sec-eval}
Our aim in this section is to prove the following result:
\begin{Thm}\label{thm-eval}
  Let $\nu\in\cVqm$ be a quasimonomial valuation. 
  There exists a unique real-valued function
  $T_\nu$ on the set of psh functions, with the following 
  properties:
  \begin{itemize}
  \item
    \emph{Compatibility:} $T_\nu(\log|\psi|)=\nu(\psi)$ for $\psi\in R$;
  \item
    \emph{Monotonicity:} $T_\nu(u)\le T_\nu(v)$ if $u\ge v+O(1)$;
  \item
    \emph{Homogeneity:}
    $T_\nu(su)=s\,T_\nu(u)$ for $s\ge0$;
  \item
    \emph{Tropicality:}
    $T_\nu(u+v)=T_\nu(u)+T_\nu(v)$;
    $T_\nu\max\{u,v\}=\min\{T_\nu(u),T_\nu(v)\}$;
  \item
    \emph{Semi-continuity:} $\limsup T_\nu(u_n)\le T_\nu(u)$  
    when $u_n\to u$ in $\Lloc$;
  \item
    \emph{Minimality:} if $T'_\nu$ satisfies all the 
    properties above, then $T_\nu\le T'_\nu$.
  \end{itemize}
  We shall then simply write $T_\nu(u)=\nu(u)$.
\end{Thm}
\begin{Remark}
  The term ``tropicality'' refers to the fact that the two 
  spaces consisting of psh functions, and nonnegative real numbers,
  both carry a tropical semi-ring structure:
  ``multiplication'' is given by addition and
  ``addition'' by $\max/\min$. 
  Thus $T_\nu$ is a homomorphism of semi-rings.
\end{Remark}
\begin{Remark}
Note that $\nu$ may be evaluated on currents: if $S=dd^cu$, then we
define $\nu(S)=\nu(u)$.  By tropicality and monotonicity, this does
not depend of the choice of $u$. 
\end{Remark}
\begin{Remark}
When a valuation $\nu$ is equivalent to a normalized valuation in the
valuative tree, \ie $\nu=c\,\nu'$ with $c>0$ and $\nu'\in\cV$, then we
set $\nu(u)\=c\,\nu'(u)$.  All properties listed above remain valid.
\end{Remark}
Our approach is as follows. We give two definitions of $T_\nu(u)$:
first as the growth rate of $u$ in a characteristic region, then 
as a generalized Lelong number. Both definitions depend on a 
choice of representation $\nu=\nu_{\phi,t}$ and a smooth germ $x$
transverse to $\phi$. Nevertheless, we show that they give the 
right value on psh functions with logarithmic singularities.
Using Demailly approximation we then show that the definitions
of $T_\nu$ in fact do not depend on the choice of $\phi$ and $x$.

Throughout the analysis, we will distinguish between the 
\emph{monomial} case $m=m(\nu)=1$ and the \emph{nonmonomial} case $m>1$. 
In the former, $T_\nu$ can be viewed as a Kiselman number and 
Theorem~\ref{thm-eval} translates into essentially well-known properties 
of Kiselman numbers. As we show, the case $m>1$ can then be
reduced to the case $m=1$ through a carefully devised 
\emph{monomialization} procedure.

\smallskip
For the rest of this section we fix a quasimonomial 
valuation $\nu\in\cVqm$.
%
%%%%%%%%%%%%%%%%%%%%%%%%%%%%%%%%%%%%%%%%%%%%%%%%%%%%%%%%%%%%%%%%%%%%%%%%%%%
%
\subsection{Characteristic regions}\label{def-qm}
Write $\nu=\nu_{\phi,t}$ for some irreducible $\phi\in\fm$. 
Unlike the skewness $t$, the germ $\phi$ 
is not uniquely determined by $\nu$, 
but we may---and will---assume that $\phi$ has minimal
multiplicity, \ie $m(\phi)=m(\nu)=:m$.
Pick $x\in\fm$ with $m(x)=1$ and $\nu_x\cdot\nu_\phi=1$,
\ie $\{x=0\}$ is smooth and transverse to $\{\phi=0\}$.
\begin{Def}
  An open set of the form
  \begin{equation}\label{e901}
    \Omega_{\phi,t,x}(r):=\{|x|<r,\ |\phi|<|x|^{mt}\},
  \end{equation}
  for small $r>0$, 
  is called a \emph{characteristic region} for $\nu$:
  see Figure~\ref{F3} on p.\pageref{F3}.
\end{Def}
Our goal is to show:
\begin{Prop}\label{P2}
  If $u$ is psh, then the limit 
  \begin{equation}\label{e2}
    \nu_{\phi,t,x}(u)
    :=\lim_{r\to0}\frac1{\log r}\sup_{\Omega_{\phi,t,x}(r)}u
  \end{equation}
  exists. Moreover
  \begin{equation}\label{e4}
    u(q)\le\nu_{\phi,t,x}(u)\log\|q\|+O(1)
  \end{equation}
  for all $q$ in $\Omega_{\phi,t,x}(r)$.
\end{Prop}
The dependence of the characteristic region on the choices 
of $\phi$ and $x$ is quite weak. In fact, we will show later
that the quantity $\nu_{\phi,t,x}(u)$ does not depend on
these choices. 
In doing so, it will be important to control the volumes of the 
characteristic regions:
\begin{Prop}\label{P3}
  There exists $C=C(\phi,t,x)>0$ such that
  \begin{equation*}
    C^{-1}r^{2A}\le\vol\Omega_{\phi,t,x}(r)\le Cr^{2A}
  \end{equation*}
  for small $r$, where $A=A(\nu)$ is the thinness of valuation 
  $\nu=\nu_{\phi,t}$.
\end{Prop}
For the proofs, we consider the cases $m=1$ and $m>1$ separately.
%
%%%%%%%%%%%%%%%%%%%%%%%%%%%%%%%%%%%%%%%%%%%%%%%%%%%%%%%%%%%%%%%%%%%%%%%%%%%
%
\subsection{The monomial case}\label{def-mono}
We first assume that $m(\nu)=1$. Then $(x,\phi)$ define local
coordinates. Write $y=\phi$ for definiteness. In 
coordinates $(x,y)$, the valuation $\nu=\nu_{y,t}$ is monomial
and the characteristic region is of the form
\begin{equation*}
  \Omega_{\phi,t,x}(r)=\{|x|<r,\ |y|<|x|^t\},
\end{equation*}
the volume of which is given by $Cr^{2+2t}$ for some constant
$C>0$. This proves Proposition~\ref{P3} in this case as $A(\nu)=1+t$.

Further, the existence of the limit in~\eqref{e2} and the 
inequality~\eqref{e4} both follow from the discussion of Kiselman
numbers in Section~\ref{sec-kiselman}.
Indeed, we have 
\begin{equation}\label{e907}
  \nu_{\phi,t,x}(u)=\nu^{x,y}_{1,t}(u)
\end{equation}
in this case, \ie $\nu$ is a Kiselman
number in coordinates $(x,y)$.

%
%%%%%%%%%%%%%%%%%%%%%%%%%%%%%%%%%%%%%%%%%%%%%%%%%%%%%%%%%%%%%%%%%
%
\subsection{Monomialization}\label{sec-precise-mono}
Many statements about quasimonomial valuations can be reduced to
the monomial case through a procedure that we will refer to
as \emph{monomialization}. Since this procedure is of fundamental
importance to our analysis, we discuss it in some detail.

Consider a quasimonomial valuation $\nu$ and write $\nu=\nu_{\phi,t}$
with $m(\phi)=m(\nu)=m$ as above. Also pick a transverse germ $x$.
We will assume that $m>1$. This implies $t>1$.
Consider the approximating sequence
$\nu_\fm=\nu_0<\nu_1<\dots<\nu_g<\nu$ of $\nu$ as in~\eqref{e705}.
The valuation $\nu_g$ is divisorial and of the form $\nu_{\phi,t_0}$ 
for some $t_0\in(1,t)$. It has generic multiplicity $m$.

Let $\pi:X\to(\C^2,0)$ be the minimal desingularization of the curve
$C=\{\phi=0\}$. This can be constructed as follows.  Let
$p_1,p_2,\dots$ be the sequence of infinitely nearby points associated
to $C$. They are defined recursively: $p_{j+1}$ is the intersection of
the strict transform of $C$ with the exceptional divisor of the blowup
at $p_j$.  We denote by $p=p_{n+1}$ the first point for which the
strict transform of $\phi$ is smooth and transverse to the exceptional
divisor at $p$.  In suitable coordinates this strict transform is
given by $\{w=0\}$, and the exceptional divisor at $p$ by $E=\{z=0\}$.

\begin{figure}[ht]
  \includegraphics[width=0.85\textwidth]{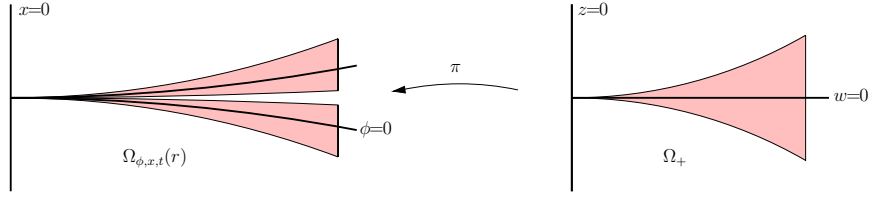}
  \caption{Monomialization. See Proposition~\ref{blw-jac}}\label{F6}
\end{figure}

Then $\pi$ is the composition of blowups at 
$p_0,\dots,p_n$. Write $J\pi$ for the Jacobian determinant of $\pi$.  
For $s\in[0,\infty)$, write $\mu_s$ for the monomial valuation sending $z$
to $1$ and $w$ to $s$. Notice that $\mu_0=\div_E$.  
By~\cite[Corollary~6.42]{treeval} the divisorial valuation
associated to $E$ is precisely $\nu_g$.
In particular, $\pi_*\div_E=m\nu_g$.
\begin{Prop}\label{blw-jac}
  Let $m$, $t_0$ and $A_0$ be the generic multiplicity, skewness
  and thinness of $\nu_g$, respectively. Then the following hold:
  \begin{itemize}
  \item[(i)]
    if $s\ge0$ then $\pi_*\mu_s=m\nu_{\phi,t_0+m^{-2}s}$;
    in particular $\pi_*\mu_{m^2(t-t_0)}=m\nu$;
  \item[(ii)]
    after multiplying $z$ and $w$ by units, if necessary, 
    we have
    \begin{equation}\label{e706}
      \pi^*x=z^m,
      \quad
      \pi^*\phi=z^{m^2t_0}w
      \qand
      J\pi=z^{mA_0-1}\xi
    \end{equation}
    for some unit $\xi$;
  \item[(iii)]
    for any $r>0$ small enough, the contraction map $\pi$ induces
    a biholomorphism from the open set
    $\Omega_+=\{(z,w),\ |z|<r^{1/m},\ |w|<|z|^{m^2(t-t_0)}\}$ 
    onto the characteristic region 
    $\Omega_{\phi,t,x}(r)$ defined in~\eqref{e901}.
    See Figure~\ref{F6}.
  \end{itemize}
\end{Prop}
The contraction map $\pi$ can also be used to show that
$\cV$ has a self-similar structure: 
see~\cite[Theorem~6.51, Figure~6.12]{treeval}.
\begin{proof}
  We start by proving~(ii). Since $x$ is transverse to $\phi$,
  the strict transform of $x$ under $\pi$ is invertible at $p$.
  As $\pi_*\mu_0=m\nu_g$ and $\nu_g(x)=1$, this gives $\pi^*x=z^m\eta_1$ 
  for some unit $\eta_1$. We may assume $\eta_1=1$ after 
  multiplying $z$ by a unit.
  Similarly, since $\nu_g(\phi)=mt_0$ and the strict transform of 
  $\{\phi=0\}$ is $\{w=0\}$, we get $\pi^*\phi=z^{m^2t_0}w\eta_2$ 
  for some unit $\eta_2$.
  By adjusting $w$ we may assume $\eta_2=1$.
  Finally, the formula for $J\pi$ follows from
  the discussion in Section~\ref{divisorial}.

  As for~(i), write $\nu_s=m^{-1}\pi_*\mu_s$ for $s\ge0$. 
  Then $\nu_s(x)=m^{-1}\mu_s(z^m)=1$ and
  $\nu_s(\phi)=m^{-1}\mu_s(z^{m^2t_0}w)=mt_0+m^{-1}s>m$ by~(ii).
  In particular, since $\{x=0\}$ and $\{\phi=0\}$ are transverse, 
  $\nu_s$ is normalized, \ie $\nu_s(\fm)=1$. 
  Since $\pi$ is a morphism,
  the assignment $s\mapsto\nu_s$ is therefore a continuous map
  of $[0,\infty]$ into $\cV$. If $s$ is irrational, then
  $\nu_s(\phi)=mt_0+m^{-1}s$ is also irrational and this 
  implies that $\nu_s=\nu_{\phi,t_0+m^{-2}s}$: 
  see Section~\ref{sec-skewness}. 
  By continuity, this must hold for all $s$, proving~(i).

  Let us finally prove the analytic statement~(iii). 
  Fix $r_0$ such that $\pi$ is well-defined and $\xi$
  invertible in the region $\Omega_+$ for $r\le r_0$
  and consider $r<r_0$.
  It follows from~(ii) that $\pi$ is an injective
  map of $\Omega_+$ into $\Omega_{\phi,t,x}(r)$. 
  Let us prove that it is surjective. 

  Fix $x$ with $|x|<r$ and let $A_x$ be the set 
  $\{y\ ;\ |\phi(x,y)|<|x|^{mt}\}$. This is a possibly 
  disconnected open subset of $\C$. 
  Let $\omega$ be one of its components.
  The proof is complete if we can
  find $z$ with $z^m=x$ such that $\pi_z$, defined by
  $\pi(z,w)=(z^m,\pi_z(w))$, maps the disk
  $\Delta=\{w\ ;\ |w|<r^{1/m}\}$ onto $\omega$. 

  Notice that $y\to |\phi(x,y)|^2$, being harmonic and non constant,
  is an open mapping, hence $\partial\omega$ is equal to $\{y\ ;\
  |\phi(x,y)|=|x|^{mt}\}$.  Moreover, $\omega$ contains
  a point $y_0$ (possibly several) such that $\phi(x,y_0)=0$. Since
  $\{w=0\}$ is the strict transform of $\{\phi=0\}$ 
  under $\pi$, there exists $z$ such that $\pi(z,0)=(x,y_0)$.
  Then $\pi_z(0)=y_0$ so since
  $\Delta$ is connected we have $\pi_z(\Delta)\subset\omega$.
  As the map $\pi_z$ is open, its image $\pi_z(\Delta)$ is open in
  $\omega$. Pick $q_n\in\pi_z(\Delta)$ converging to $q\in\omega$, and
  write $q_n=\pi_z(p_n)$, $p_n\in\Delta$. 
  Extract a subsequence $p_n\to p\in\overline{\Delta}$.  
  By continuity $\pi_z(p)=q$. But~(ii) shows
  that $\pi_z(\partial\Delta)\subset\partial\omega$, hence
  $p\in\Delta$, and $\pi_z(\Delta)$ is closed in $\omega$.  The open
  set $\omega$ being connected we conclude that
  $\pi_z(\Delta)=\omega$, which completes the proof.
\end{proof}
We are now ready to prove the existence of $\nu_{\phi,t,x}(u)$
and the volume estimate in Proposition~\ref{P3} in the
nonmonomial case $m>1$.
\begin{proof}[Proof of Proposition~\ref{P3}]
  Write $s=m^2(t-t_0)$. 
  By Proposition~\ref{blw-jac} we have
  $|J\pi(z,w)|\simeq|z|^{mA_0-1}$
  in $\Omega_+$ if $r$ is small enough.
  By the change of variables formula, this yields
  \begin{equation*}
    \vol\Omega_{\phi,t,x}(r)
    =\int_{\Omega_+}|J\pi|^2
    \simeq\int_{|z|<r^{\frac1m},|w|<|z|^s}|z|^{2mA_0-2}
    \simeq r^{m^{-1}(2+2s+2mA_0-2)},
  \end{equation*}
  which completes the proof, since 
  $A(\nu)=A_0+m(t-t_0)=A_0+m^{-1}s$.
\end{proof}
\begin{proof}[Proof of Proposition~\ref{P2}]
  Again we use Proposition~\ref{blw-jac}, which gives
  \begin{equation*}
    \frac1{\log r}\sup_{\Omega_{\phi,t,x}(r)}u
    =\frac1{\log r}\sup\{\pi^*u\ ;\ |z|<r^{\frac1m},|w|<|z|^s\}.
  \end{equation*}
  for small $r$.
  Hence $\nu_{\phi,t,x}(u)$ is well-defined and we have
  \begin{equation}\label{e903}
    \nu_{\phi,t,x}(u)=m^{-1}\nu^{z,w}_{1,s}(\pi^*u),
  \end{equation}
  where $\nu^{z,w}$ denotes the Kiselman number in coordinates $(z,w)$.
  Notice that~\eqref{e903} exhibits $\nu_{\phi,t,x}$ as the pushforward
  of a Kiselman number.
  
  As for~\eqref{e4}, pick $q=(x,y)\in\Omega_{\phi,t,x}(r)$
  and write $q=\pi(q')$ where $q'=(z,w)$.
  Let us first assume that $s\ge1$.
  Then $\|q'\|\simeq|z|=|x|^{1/m}\simeq\|q\|^{1/m}$ 
  so~\eqref{e902} implies
  \begin{equation*}
    u(q)
    =(\pi^*u)(q')
    \le\nu^{z,w}_{1,s}(\pi^*u)\log\|q'\|+O(1)
    =\nu_{\phi,t,x}(u)\log\|q\|+O(1).
  \end{equation*}
  If instead $s<1$, then
  $\|q'\|\simeq|w|=|\phi|/|z|^{m^2t_0}\le|x|^{s/m}\simeq|q|^{s/m}$, so
  \begin{align*}
    u(q)
    &\le\frac1s\nu^{z,w}_{1,s}(\pi^*u)\log\|q'\|+O(1)
    \le\frac{m}{s}\nu_{\phi,t,x}(u)\frac{s}{m}\log\|q\|+O(1)\\
    &=\nu_{\phi,t,x}(u)\log\|q\|+O(1),
  \end{align*}
  which concludes the proof.
\end{proof}
%
%%%%%%%%%%%%%%%%%%%%%%%%%%%%%%%%%%%%%%%%%%%%%%%%%%%%%%%%%%%%%%%%%
%
\subsection{Analytic definition}
We have defined $\nu_{\phi,t,x}$ as a growth rate in the
characteristic region. Let us give a more analytic definition.
\begin{Prop}\label{P6}
  Given $\phi$, $t$ and $x$ as above, we have
  \begin{equation}\label{e601}
    \nu_{\phi,t,x}(u)
    =dd^cu\wedge dd^c\log\max\{|x|^t,|\phi|^{1/m}\}\,\{0\}
  \end{equation}
  for any psh function $u$.
  Thus $\nu_{\phi,t,x}(u)$ can be viewed as the generalized 
  Lelong number $\nu_\varphi(u)$ with weight 
  $\varphi=\log\max\{|x|^t,|\phi|^{1/m}\}$.
\end{Prop}
\begin{proof}
  Write $\nu^{\mathrm{an}}_{\phi,t,x}(u)$ 
  for the right hand side of~\eqref{e601}.
  
  In the monomial case $m=m(\nu)=1$,~\eqref{e601} 
  reduces to a well-known statement about Kiselman numbers 
  by putting $y=\phi$: see~\eqref{e907} and~\eqref{e906}.

  When $m>1$ then monomialization as in Proposition~\ref{blw-jac} gives:
  \begin{align*}
    \nu^{\mathrm{an}}_{t,\phi,x}(u)
    &=\lim_{r\to0}\frac{1}{\log r}
    \int u\,(dd^c)^2\log\max\{r,|x|^t,|\phi|^{1/m}\}\\
    &=\lim_{r\to0}\frac{1}{mt\log r}
    \int(\pi^*u)\,(dd^c)^2\log\max\{r^{mt},|z|^{mt},|w|^{1/m}|z|^{mt_0}\}\\
    &=\lim_{r\to0}\frac{1}{m\log r}
    \int(\pi^*u)\,(dd^c)^2\log\max\{r,|z|,|w|^s\}\\
    &=m^{-1}\nu^{z,w}_{1,s}(\pi^*u),
  \end{align*}
  where $s=m^2(t-t_0)$.
  Here the first equality is a consequence of~\eqref{DL-jensen}, the 
  second follows from the change of variables formula and from replacing
  $r$ by $r^{mt}$. In the third equality we have used
  Lemma~\ref{L300}. The last line follows from~\eqref{e906}.
  Thus~\eqref{e903} implies that 
  $\nu^{\mathrm{an}}_{t,\phi,x}(u)=\nu_{\phi,t,x}(u)$.
\end{proof}
%
%%%%%%%%%%%%%%%%%%%%%%%%%%%%%%%%%%%%%%%%%%%%%%%%%%%%%%%%%%%%%%%%%
%
\subsection{Monotonicity, tropicality and homogeneity}\label{S102}
It is immediate from~\eqref{e2} that if $u$ and $v$ are
psh functions with $u\ge v+O(1)$, then 
$\nu_{\phi,t,x}(u)\le\nu_{\phi,t,x}(v)$.
Moreover, we have
\begin{Prop}\label{P233}
  Given $\phi$, $t$ and $x$ as above, we have
  \begin{equation*}
    \nu(su+tv)=s\,\nu(u)+t\,\nu(v)
    \qand
    \nu(\max\{u,v\})=\min\{\nu(u),\nu(v)\},
  \end{equation*}
  for any psh functions $u$,$v$ and any constants $s,t\ge0$.
  Here $\nu=\nu_{\phi,t,x}$.
\end{Prop}
\begin{proof}
  The first equality is an immediate consequence of
  Proposition~\ref{P6}.
  As $u\le\max\{u,v\}$,~\eqref{e2} 
  immediately implies $\nu(\max\{u,v\})\le\nu(u)$. 
  By symmetry we obtain 
  $\nu(\max\{u,v\})\le\min\{\nu(u),\nu(v)\}$.  
  For the reverse inequality we use~\eqref{e4}:
  in $\Omega_{\phi,t,x}(r)$ we have 
  $u(q)\le\nu(u)\log\|q\|+O(1)$ 
  and $v(q)\le\nu(v)\log\|q\|+O(1)$.
  Thus $\max\{u,v\}(q)\le\min\{\nu(u),\nu(v)\}\log\|q\|+O(1)$,
  so~\eqref{e2} immediately gives
  $\nu(\max\{u,v\})\ge\min\{\nu(u),\nu(v)\}$. 
\end{proof}
%
%%%%%%%%%%%%%%%%%%%%%%%%%%%%%%%%%%%%%%%%%%%%%%%%%%%%%%%%%%%%%%%%%
%
\subsection{Compatibility}
Next we show that independently of any choice of $\phi$ and $x$,
the definition of $\nu(u)$ agrees with the algebraic definition
when $u$ has logarithmic singularities, in the following sense.
\begin{Prop}\label{P4}
  Let $u=c\log(\sum_1^n|\psi_i|)+O(1)$, $c>0$, $\psi_i\in R$, 
  be a psh function with logarithmic singularities. Then
  \begin{equation}\label{eq-intsup1}
    \nu_{\phi,t,x}(u)=c\min_i\nu(\psi_i).
  \end{equation}
  Here $\nu(\psi_i)$ denotes the value of $\nu$ 
  on the holomorphic germ $\psi_i\in R$.
\end{Prop}
In particular, the definition of $\nu(u)$ is independent of 
the choices of $\phi$ and $x$ when $u$ has logarithmic singularities.
\begin{proof}
  We may assume $c=1$ and that the $O(1)$ term is absent.
  Moreover,
  $\min\log|\psi_i|\le\log\sum_1^n|\psi_i|\le\min\log|\psi_i|+\log n$, 
  so $\nu(\log\sum|\psi_i|)=\min\nu(\log|\psi_i|)$ 
  as follows from~\eqref{e4}. 
  We can hence suppose $u=\log|\psi|$ for $\psi\in R$.

  First consider the monomial case $m=m(\nu)=1$. Write $y=\phi$
  for definiteness. 
  Then $\psi$ can be viewed as a power series in
  coordinates $(x,y)$, say
  $\psi=\sum a_{kl}x^ky^l$.
  By definition $\nu(\psi)=\min\{ k+tl\ ;\ a_{kl}\ne0\}=:\a$.
  The characteristic region is given by
  $\Omega_{\phi,t,x}(r)=\{|x|<r,|y|<r^t\}$ and it is elementary
  to see that $\sup |\psi|\simeq r^\a$ on $\Omega_{\phi,t,x}(r)$.
  Therefore,
  $\nu_{\phi,t,x}(\log|\psi|)=\a=\nu(\psi)$.
  
  When $m>1$, we monomialize, using 
  Proposition~\ref{blw-jac}. Write $s=m^2(t-t_0)$.
  Then $\nu=m^{-1}\pi_*\mu_s$ so by~\eqref{e903} and
  the monomial case above, we get
  \begin{equation*}
    \nu_{\phi,t,x}(\log|\psi|)
    =m^{-1}\nu^{z,w}_{1,s}(\log|\pi^*\psi|)
    =m^{-1}\mu_s(\pi^*\psi)
    =\nu(\psi).
  \end{equation*}
  This concludes the proof.
\end{proof}
%
%%%%%%%%%%%%%%%%%%%%%%%%%%%%%%%%%%%%%%%%%%%%%%%%%%%%%%%%%%%%%%%%%
%
%
\subsection{Approximation}
To show that the definition of $\nu(u)$, for a general psh function $u$,
is independent of all choices made, we use Demailly approximation.
\begin{Prop}\label{P5}
  For any psh function $u$, there exists a sequence $u_n$ of psh
  functions with logarithmic singularities such that $u_n\to u$ in
  $\Lloc$ and 
  \begin{equation}\label{e904}
    0\le\nu_{\phi,t,x}(u)-\nu_{\phi,t,x}(u_n)\le\frac{A}{n}
  \end{equation}
  for any choice of $\phi$ and $x$. Here $A=A(\nu)$
  is the thinness of $\nu=\nu_{\phi,t}$.
\end{Prop}
\begin{Remark}
  The estimate~\eqref{e904} for Lelong numbers (\ie $\nu=\nu_\fm$) is
  due to Demailly (see~\cite[Theorem~4.2]{DK}). The extension to
  Kiselman numbers (\ie $m(\nu)=1$) was proved independently by
  Rashkovskii~\cite[Theorem~3]{rash} and the
  authors~\cite[Lemma~2.4]{FJ-brolin}.
\end{Remark}
\begin{proof}
  We choose $u_n$ as in Section~\ref{sec-demailly},
  \ie $u_n=n^{-1}\sup\log|f|$ where the supremum is over
  holomorphic functions $f$ such that $\int|f|^2\exp(-2nu)\le1$ 
  over a fixed neighborhood of the origin.
  
  By~\eqref{e801}, $u_n\ge u-C$ for some constant $C>0$, hence
  $\nu_{\phi,t,x}(u_n)\le\nu_{\phi,t,x}(u)$. For the other inequality,
  consider the open set
  \begin{equation*}
    \Omega'(r)=\Omega'_{\phi,t,x}(r)
    =\{\frac12r<|x|<r,\ \frac12|x|^{mt}<|\phi|<|x|^{mt}\}.
  \end{equation*}
  With minor modifications, the proof of Proposition~\ref{P3} shows that 
  $\vol\Omega'(r)\simeq r^{2A(\nu)}$ and by the maximum principle
  $(\log r)^{-1}\sup_{\Omega'(r)}v=(\log r)^{-1}\sup_{\Omega(r)}v
  \to\nu_{\phi,t,x}(v)$ 
  as $r\to0$, for any psh function $v$.
  We need the following estimate:
  \begin{Lemma}\label{L901}
    For any $\e\in(0,1)$ there exists $C_\e>1$ such that
    \begin{equation*}
      C_\e\sup_{\Omega'(\e r)}v
      \le
      \frac{1}{\vol\Omega'(r)}\int_{\Omega'(r)}v
    \end{equation*}
    for small $r$ and any psh function $v\le0$. Moreover, 
    $C_\e\to1$ as $\e\to0$.
  \end{Lemma}
  We apply this lemma to $v=\log|f|^2$, where 
  $\int|f|^2\exp(-2nu)\le1$.
  Using the concavity of the logarithm we get
  \begin{multline*}
    C_\e\sup_{\Omega'(\e r)}\log|f|^2
    \le\frac{1}{\vol\Omega'(r)}\int_{\Omega'(r)}\log|f|^2
    \le\log\left[\frac{1}{\vol\Omega'(r)}\int_{\Omega'(r)}|f|^2\right]\le\\
    \le\log\left[
      \frac{1}{\vol\Omega'(r)}\sup_{\Omega'(r)}e^{2nu}
      \int_{\Omega'(r)}|f|^2e^{-2nu}
    \right]
    =2n\sup_{\Omega'(r)}u-\log\vol\Omega'(r).
  \end{multline*}
  As $\log\vol\Omega'(r)=2A(\nu)\log r+O(1)$ 
  we get $C_\e\nu_{\phi,t,x}(\log|f|)\ge n\nu_{\phi,t,x}(u)-A(\nu)$
  by letting $r\to0$.
  Letting $\e\to0$ yields 
  $\nu_{\phi,t,x}(\log|f|)\ge n\nu_{\phi,t,x}(u)-A(\nu)$.
  Using $u_n=n^{-1}\sup\log|f|$, we conclude 
  $\nu_{\phi,t,x}(u_n)\ge\nu_{\phi,t,x}(u)-A(\nu)/n$.
\end{proof}
\begin{proof}[Proof of Lemma~\ref{L901}]
  We only prove this result in the case $m>1$, the
  monomial one ($m=1$) being easier.  
  As usual we make use of the monomialization procedure in
  Proposition~\ref{blw-jac}. Write
  \begin{equation*}
    \cD(r)
    =\{(\frac12r)^{\frac1m}<|z|<r^{\frac1m},\ \frac12|z|^{s}<|w|<|z|^s\}.
  \end{equation*}
  Then $\pi$ is a biholomorphism from $\cD(r)$ onto $\Omega'(r)$.
  Let $D(r)=\{(|z|,|w|)\ ;\  (z,w)\in\cD(r)\}\subset\R^2$.
  For $\rho=(\rho_1,\rho_2)\in D(r)$ let $T_\rho$ be the torus
  $\{|z|=\rho_1,|w|=\rho_2\}\subset\cD(r)$.
  For any $\e\in(0,1)$ and any $\rho\in D(r)$, 
  Lemma~\ref{L-harnack} and the maximum principle imply
  \begin{equation*}
    C_\e\sup_{\Omega'(\e r)}v
    =C_\e\sup_{\cD(\e r)}\pi^*v
    \le\int_{T_\rho}\pi^*v,
  \end{equation*}
  where the integral is with respect to Lebesgue measure
  of mass one on the torus.
  Note that we can take $C_\e = (1+\e)^4/ (1-\e^2)$ 
  hence  $C_\e\to 1$ when $\e\to0$.
  
  The critical set of $\pi$ is equal to $\{z=0\}$.
  We can thus write $C'|z|^{mA_0-1}\ge |J\pi|\ge C|z|^{mA_0-1}$  on $\cD(r)$,
   for some constants $C = C(r),C'= C'(r)>0$, with
   $\lim_{r\to0}C(r)/C'(r)=1$. We have
  \begin{multline*}
    \int_{\Omega'(r)} v
    =\int_{\cD(r)} |J\pi|^2\pi^*v
    =\int_{\rho\in D(r)} 
    \left[
      \int_{T_\rho} |J\pi|^2\pi^*v
    \right]
    \rho_1\rho_2\,d\rho_1d\rho_2\\
    \ge\int_{\rho\in D(r)} 
    \left[
      \int_{T_\rho} C^2\rho_1 ^{2(mA_0-1)}\pi^*v
    \right]
    \rho_1\rho_2\,d\rho_1d\rho_2\\
    \ge (C/C')^2 C_\e\int_{\rho\in D(r)} 
    \left[
      \int_{T_\rho} |J\pi|^2 
      \sup_{\cD(\e r)}\pi^*v 
    \right]
    \rho_1\rho_2\,d\rho_1d\rho_2
    =C'_\e\vol\Omega'(r)\,\sup_{\Omega'(\e r)}v,
  \end{multline*}
  with $C'_\e\to 1$ as $\e\to0$.
\end{proof}
%
%%%%%%%%%%%%%%%%%%%%%%%%%%%%%%%%%%%%%%%%%%%%%%%%%%%%%%%%%%%%%%%%%
%
\subsection{Proof of Theorem~\ref{thm-eval}}
Let us summarize our construction.  Fix a quasimonomial valuation $\nu$. 
Write $\nu=\nu_{\phi,t}$ with $m(\phi)=m(\nu)=:m$ and pick
a holomorphic germ $x$ such that $\{x=0\}$ is smooth and
transverse to $\{\phi=0\}$. We defined $T_\nu$
in two equivalent ways: as a growth rate on a characteristic
region (Proposition~\ref{P2}) and as a generalized Lelong
number (Proposition~\ref{P6}). These definitions a priori
depended on the choices of $\phi$ and $x$ but we showed independence
of these choices when restricted to psh functions with logarithmic
singularities (Proposition~\ref{P4}). Using Demailly approximation
(Proposition~\ref{P5}), this was extended to arbitrary psh functions.
Monotonicity, homogeneity and tropicality were proved in
Section~\ref{S102}. Semi-continuity is an immediate consequence 
of the analytic definition of $T_\nu$. Compatibility was proved
in Proposition~\ref{P4}. 
Finally, let us prove minimality (and hence uniqueness).
By compatibility we have $T'_\nu(u)=T_\nu(u)$ when $u=\log|\phi|$,
$\phi\in R$. By Proposition~\ref{P4}, tropicality, 
and monotonicity, the same result holds when $u$ has logarithmic
singularities. Semi-continuity together with Proposition~\ref{P5}
then shows that $T'_\nu(u)\ge T_\nu(u)$ for all $u$.
This concludes the proof of Theorem~\ref{thm-eval}.
%
%
%%%%%%%%%%%%%%%%%%%%%%%%%%%%%%%%%%%%%%%%%%%%%%%%%%%%%%%%%%%%%%%%%
%
%
\section{Further properties}\label{sec-further}
In this section we give further information on the
pairing $(\nu,u)\mapsto\nu(u)$.
%
%%%%%%%%%%%%%%%%%%%%%%%%%%%%%%%%%%%%%%%%%%%%%%%%%%%%%%%%%%%%%%%%%
%
\subsection{Geometric interpretation}
Let us prove the following geometric interpretation of $\nu(u)$ in
case $\nu$ is divisorial.
\begin{Prop}\label{P-interpret}
  Suppose $\nu\in\cV$ is a divisorial valuation, associated to an
  exceptional component $E$ of a modification $\pi$ above the
  origin. Write $\nu=\nu_E=b^{-1}\pi_*\div_E$, with 
  $b$ the generic multiplicity of $\nu$.
  
  Then for any psh function $u$, $\nu(u)$ is equal to $1/b$ times the
  Lelong number of $\pi^*u:=u\circ\pi$ taken at a generic point of $E$.
\end{Prop}
\begin{proof}
  By Siu's theorem, the Lelong number of $\pi^*u$ at all points 
  $p\in E$ except for countably many is the same, 
  equal to some constant $M\ge0$.
  Pick a point $p\in E$ in the regular set of $\pi^{-1}(0)$ such that
  $\nuL(\pi^*u,p)=M$.  Choose a smooth curve $V$ cutting $E$
  transversely at $p$. We let $\phi\in\fm$ be an equation defining
  $\pi(V)$. By Section~\ref{divisorial}, $\nu_\phi\ge\nu$ hence
  $\nu=\nu_{\phi,t}$ for some $t$.  
  Fix local coordinates $(z,w)$ at $p$, such that 
  $E=\{z=0\}$, and $V=\{w=0\}$, and a coordinate axis $x$
  at the origin in $\C^2$ such that Proposition~\ref{blw-jac}~(iii)
  applies: the map $\pi$ gives a biholomorphism between 
  $\Delta_{r^{1/b}}\times\Delta_1$ and 
  $\Omega_{\phi,x,t}(r)=\{|x|<r,|\phi|<|x|^{bt}\}$ 
  (note that in Proposition~\ref{blw-jac} the
  generic multiplicity of $\nu$ is denoted by $m$).
  We hence obtain
  $\sup_{\Omega_{\phi,t,x}(r)}u=\sup_{\Delta_{r^{1/b}}\times\Delta_1}\pi^*u$
  and
  \begin{equation*}
    \nu(u)
    =\lim_{r\to0}\frac1{\log r}
    \sup_{\Delta_{r^{1/b}}\times\Delta_1}\pi^*u.
  \end{equation*}
  As 
  $\Delta_{r^{1/b}}\times\Delta_{r^{1/b}}
  \subset\Delta_{r^{1/b}}\times\Delta_1$
  we get $b\times\nu(u)\le\nuL(\pi^*u,p)=M$.
  
  It is also a consequence of Siu's theorem that the function 
  $\pi^*u-M\log|z|$ is psh. By the maximum principle, the supremum of 
  $\pi^*u-M \log|z|$ on 
  $\Delta_{r^{1/b}}\times\Delta_1$ is attained on 
  $\partial\Delta_{r^{1/b}}\times\Delta_1$, hence
  \begin{multline*}
    0\le\lim_{r\to0} \frac1{\log r}
    \sup\left\{\pi^*u-M \log |z| \ ; \ |z|<  r^{1/b} \ ,
      |w|<1 \right\}
    =\\
    =\lim_{r\to0} \frac1{\log r}
    \sup\left\{\pi^*u-\frac{M}{b}\log r\ ;\ |z|=r^{1/b},|w|<1\right\}
    = \\
    =-\frac{M}{b}
    +\lim_{r\to0}\frac1{\log r}
    \sup\left\{\pi^*u\ ; |z|<r^{1/b},|w|<1\right\}
    =-\frac{M}{b}+\nu(u)
    \le 0
  \end{multline*}
We conclude that $b\times\nu(u)=M$.
\end{proof}
%
%%%%%%%%%%%%%%%%%%%%%%%%%%%%%%%%%%%%%%%%%%%%%%%%%%%%%%%%%%%%%%%%%
%
\subsection{Action by holomorphic maps}
Next we prove that the evaluation map $(\nu,u) \to \nu(u)$ behaves
well under the action of a holomorphic map, a fact relevant to 
the study of the action of holomorphic maps on the
valuative tree~\cite{eigenval}.
\begin{Prop}
  Suppose $f : (\C^2,0) \to (\C^2,0)$ is a holomorphic germ whose
  Jacobian determinant does not vanish identically.
  Then for any plurisubharmonic function $u$ and any quasimonomial valuation
  $\nu$, we have
  \begin{equation}\label{act-fonct}
    \nu(f^*u)=(f_*\nu)(u),
  \end{equation}
  where $f^*u\=u\circ f$, and $f_*\nu$ is the valuation defined by 
  $(f_*\nu)(\phi)=\nu(f^*\phi)$ for $\phi\in R$.
\end{Prop}
Note that in general $f_*\nu$ is no longer normalized: there is
no reason why $c(\nu):=(f_*\nu)(\fm)=1$. 
It is however equivalent to a unique normalized 
valuation $f_\bullet\nu \in \cV$, and we can set 
$(f_*\nu)(u)\=c(\nu)\times(f_\bullet\nu)(u)$ for any psh function $u$.
\begin{proof}
  If $u=c\log\sum_1^n|\phi_i|^2$, where $c>0$ and $\phi_i\in R$
  are holomorphic germs, then~\eqref{act-fonct} is an immediate
  consequence of Proposition~\ref{P4}. Thus~\eqref{act-fonct}
  holds when $u$ has logarithmic singularities
  in view of Proposition~\ref{P233}.

  Now consider a general plurisubharmonic function $u$
  and approximate it by psh functions $u_n$ with 
  logarithmic singularities as in Proposition~\ref{P5}. 
  Recall that $u_n=n^{-1}\sup\log|h|$
  where $h$ ranges over all holomorphic functions such that 
  $\int_B|h|^2\exp(-2nu)\le 1$, and $B$ is a fixed ball centered 
  at $0$. Define in the same way
  $\tu_n\=n^{-1}\sup\log|\th|$ where $\th$ ranges over all holomorphic
  functions such that $\int_B|\th|^2\exp(-2nf^*u)\le1$. 
  This sequence approximates $f^*u$. By choosing suitable
  coordinates at the source space, we may assume that 
  $f(B)\subset B$. Denote by $Jf$ the Jacobian determinant of $f$. 
  The change of variables formula yields
  \begin{equation*}
    \int_B|f^*h\cdot Jf|^2\exp(-2nf^*u)
    \le e\int_{f(B)}|h|^2\exp(-2nu)
    \le e\int_{B}|h|^2\exp(-2nu),
  \end{equation*}
  where $e$ denotes the topological degree of $f$, \ie the cardinality
  of a generic fiber of $f$.  We thus infer $\tu_n\ge f^*u_n +
  n^{-1}\log|Jf/\sqrt{e}|$ for all $n$.  Using Proposition~\ref{P5} we
  conclude $\nu(f^*u)=\lim\nu(\tu_n)\le\liminf\nu(f^*u_n)$.  On the
  other hand, $f^*u_n\to f^* u$ in $\Lloc$, so by semi-continuity
  $\nu(f^*u)\ge\limsup\nu(f^*u_n)$.  This shows that
  $\nu(f^*u)=\lim\nu(f^*u_n)$, so
  \begin{equation*}
    (f_*\nu)(u)
    =\lim_{n\to\infty}(f_*\nu)(u_n)
    =\lim_{n\to\infty}\nu(f^*u_n)
    =\nu(f^*u),
  \end{equation*}
  which concludes the proof.
\end{proof}
%
%%%%%%%%%%%%%%%%%%%%%%%%%%%%%%%%%%%%%%%%%%%%%%%%%%%%%%%%%%%%%%%%%
%
\subsection{Estimates}
The characteristic regions we used to defined $\nu(u)$ can be replaced
by other semi-analytic regions.  For constants
$C_1,C_2,C_3,C_4>0$ introduce 
$\cA_{\phi,t,x,C}(r)=\{ C_1r<|x|<C_2r\ ;\ C_3|x|^{mt}<|\phi|<C_4|x|^{mt}\}$.
With these regions we can strengthen
Proposition~\ref{P2} as follows.
\begin{Prop}\label{P230}
  For any psh function $u$ and any $\nu\in\cVqm$:
  \begin{align*}
    \nu(u)
    &=\lim_{r\to0}\frac1{\log r}\sup_{\cA_{\phi,t,x,C}(r)}u
    \quad\text{and}\\
    u(q)
    &\le\nu(u)\log\|q\|+O(1)
    \quad\text{for all $q\in\cA_{\phi,t,x,C}(r)$}.
  \end{align*}
  When $u$ has logarithmic singularities, we may pick 
  $C=C(u)$ such that 
  \begin{equation*}
    u(q) 
    \ge\nu_{\phi,t,x}(u)\log\|q\|+O(1)
    \quad\text{for all $q\in\cA_{\phi,t,x,C}(r)$}.
  \end{equation*}
\end{Prop}
\begin{proof}
  The first two assertions are immediate consequences of the maximum
  principle.  For the last assertion, we may suppose $u=\log|\psi|$,
  $\psi\in\fm$ and $\nu$ is a monomial valuation sending $x$ to $1$ and
  $y$ to $t$. Write $\psi(x,y)=\sum a_{kl} x^ky^l$ and $\nu=\nu(\psi)$.
  Then $|\sum_{k+tl>\nu} a_{kl} x^ky^l|\le \mathrm{Const} \cdot
  r^{\nu+\e}$ for some $\e>0$. 
  The set of $(k,l)$ such that $k+lt=\nu$ and
  $a_{kl}\ne0$ is non-empty and finite. 
  We denote by $(k_0,l_0)$ the element of this set 
  with maximal first component (hence minimal second
  component), and we have 
  $|a_{k_0l_0}x^{k_0}y^{l_0}|\ge|a_{k_0l_0}|C_1^{k_0}C_3^{l_0}r^{\nu}$.  
  For the other $(k,l)$ with $k+lt=\nu$, we have 
  $|a_{kl}x^{k}y^{l}|\le|a_{kl}|C_2^{k}C_4^{l}r^{\nu}$. 
  Now choose $C_2>C_1\gg1$ large but
  sufficiently close to each other, and $C_3<C_4\ll1$ small.  
  With these choices of $C_i$'s, we have 
  $|\sum_{k+lt=\nu}a_{kl}x^ky^l|\ge
  |a_{k_0l_0}x^{k_0}y^{l_0}| - \sum_{k+lt=\nu, \
    (k,l)\ne(k_0,l_0)}|a_{kl}x^ky^l|\ge \mathrm{Const}\cdot
  r^{\nu}$. Therefore, $|\psi|\ge\mathrm{Const}\cdot r^{\nu}$.
  This completes the proof.
\end{proof}
%
%
%%%%%%%%%%%%%%%%%%%%%%%%%%%%%%%%%%%%%%%%%%%%%%%%%%%%%%%%%%%%%%%%%%%%%%%%%%%
%
%
\section{Potential theory on trees}\label{sec-pot}
The next step in our approach is to analyze, for a fixed 
psh function $u$, the function $\nu\mapsto\nu(u)$ on
the tree $\cVqm$ of quasimonomial valuations.

As we showed in~\cite[Chapter 7]{treeval}, there is a general
correspondence between positive measures on trees endowed with a
parameterization and certain functions called \emph{tree potentials}.
In this section, we review briefly this correspondence in the case
of the valuative tree $\cV$ parameterized by skewness, referring
to~\cite{treeval} for the proofs. Later we shall
show that the function $\nu\mapsto\nu(u)$ is a tree potential.
%
%%%%%%%%%%%%%%%%%%%%%%%%%%%%%%%%%%%%%%%%%%%%%%%%%%%%%
%
\subsection{Borel measures}\label{measure}
(See~\cite[Section~7.3]{treeval}.)  We equip $\cV$ with the Borel
$\sigma$-algebra generated by the weak topology (see
Section~\ref{S11}).  We let $\cM$ be the set of positive\footnote{Note
that in~\cite{treeval}, $\cM$ denotes the set of \emph{complex} Borel
measures.}  Borel measures\footnote{Every positive finite measure on
$\cV$ is a Radon measure by~\cite[Proposition~7.14]{treeval}.} on
$\cV$, that is, continuous positive linear functionals on the set of
continuous real valued functions on $\cV$.  The space $\cM$ also
carries a weak topology: $\rho_k\to\rho$ iff
$\int\varphi\,d\rho_k\to\int\varphi\,d\rho$ for all $\varphi$.  The
subset of probability measures is compact.  We identify $\nu\in\cV$
with the corresponding point mass $\delta_\nu\in\cM$.

A \emph{subtree} $\cT$ of $\cV$ is a subset of $\cV$ such that
$\mu\in\cT$ and $\nu\le\mu$ imply $\nu\in\cT$.  
It is a \emph{finite} subtree if it has finitely many ends.
The weak topology on $\cV$ is the weakest topology restricting 
to the usual topology on any finite subtree.

While the valuative tree has ample branching, 
the support of a measure is always much thinner:
\begin{Lemma}\label{supp-count}
  The support of any positive measure $\rho\in\cM$ is 
  contained in the closure of a countable union of finite subtrees.
\end{Lemma}
\begin{proof}
  We may suppose that $\rho$ has mass $1$.  Consider the decreasing
  function $f=f_\rho:\cV\to[0,1]$ defined by
  $f(\nu)=\rho\{\mu\ge\nu\}$. The support of $\rho$ is included in
  the closure of the union of the trees 
  $\cT_n:=\{f\ge n^{-1}\}$, $n\ge1$. 
  By construction the tree $\cT_n$ has at most $n+1$ ends.
\end{proof}
\begin{Remark}\label{R701}
  Any positive measure is uniquely 
  determined by its values on the 
  sets $U(\vv)$: see~\cite[Lemma~7.18]{treeval}.
  This is a nontrivial assertion, despite the fact that
  the weak topology is by definition generated by the sets
  $U(\vv)$.
\end{Remark}

% 
%%%%%%%%%%%%%%%%%%%%%%%%%%%%%%%%%%%%%%%%%%%%%%%%%%%%% 
%
\subsection{Tree potentials.}  
(See~\cite[Section 7.9]{treeval}).  We now describe the class of
functions on $\cVqm$ that model the behavior of tree transforms of psh
functions.

Consider a function $g:\cVqm\to\R$, and pick a tangent vector
$\vv$ at $\nu\in\cVqm$. 
We define the \emph{derivative of $g$ along $\vv$}
(when it exists) by
\begin{equation*}
  D_\vv g
  =\lim\frac{g(\nu)-g(\mu)}{|\a(\nu)-\a(\mu)|}\
  \text{when $\mu$ tends to $\nu$ along $\vv$}.
\end{equation*}
Here a sequence $\mu_k$ tends to $\nu$ along $\vv$ iff 
$\mu_k\in U(\vv)$, and the segments $[\mu_k,\nu[$ form a decreasing
sequence with empty intersection.
This is stronger than saying
that $\mu_k\to\nu$ and $\mu_k\in U(\vv)$ for all $k$.

Note that if $g$ is increasing on $\cV$, and the derivative of
$g$ along $\vv$ is well defined, then $D_\vv g \ge 0$ when
$\vv$ is not represented by $\nu_\fm$, and $D_\vv g\le0$ otherwise.
\begin{Def}\label{treepot}
  A function $g:\cVqm\to\R$ is called a \emph{tree potential} on $\cV$,
  if the following conditions are satisfied:
  \begin{itemize}
  \item[(P1)] 
    $g$ is nonnegative, increasing, and concave along totally ordered 
    segments;
  \item[(P2)]
    if $\nu\ne\nu_\fm$, then
    $\sum_{\vv\in T\nu}D_\vv g\le0$;
  \item[(P3)]
    $\sum_{\vv\in T\nu_\fm}D_\vv g\le g(\nu_\fm)$.   
  \end{itemize}  
  We denote by $\cP$ the set of all tree potentials 
  on $\cV$.\footnote{Note
    that in~\cite{treeval}, $\cP$ is denoted by $\cP^+$, and that tree
    potentials are \emph{positive tree potentials}.} 
\end{Def}
If $\nu_\star$ is an end in $\cV$, then~(P1) allows us
to define $g(\nu_\star)\in\Rbar$ as well as $D_\vv g$,
where $\vv$ is the unique tangent vector at $\nu_\star$.

As the tangent space is uncountable at any divisorial
valuation, conditions~(P2) and~(P3) are quite strong.
Indeed, if $\nu$ is divisorial, then for all but countably
many $\vv\in T\nu$ we must have $D_\vv g=0$, which 
by~(P1) implies that $g\equiv g(\nu)$
on $U(\vv)$ if $\vv$ is not represented by $\nu_\fm$.
Along these lines, we can show
\begin{Lemma}\label{lem-support}
  Suppose $g$ is a tree potential. Then its support
  \begin{equation*}
    \supp g:=\{\nu\in\cVqm\ |\ g\ \text{not locally constant at $\nu$}\},
  \end{equation*}
  is contained in the closure of a countable tree.
\end{Lemma}
\begin{proof}
  After multiplication by a constant we may
  assume that $g(\nu_\fm)\le1$.  
  For $n\ge1$, let $\cT_n$ be the set
  consisting of $\nu_\fm$ and all valuations $\nu\ne\nu_\fm$ for which
  $D_\vv g\le-n^{-1}$, where $\vv$ is the tangent vector at $\nu$
  represented by $\nu_\fm$.  Then $(\cT_n)_{n\ge1}$ forms an
  increasing sequence of trees such that the closure of their union
  equals $\supp g$. Hence we are done if we can show that $\cT_n$ is
  finite.
  
  The key remark is that we may de-localize the integrability condition
  which defines a tree potential. 
  Specifically, pick tangent vectors $\vv_1,\dots,\vv_k$ 
  at different valuations
  $\nu_1,\dots,\nu_k$, such that $\vv_i$ is represented neither by
  $\nu_\fm$, nor by $\nu_j$ for $j\ne i$. 
  Then $\sum_1^k D_{\vv_i}g\le g(\nu_\fm)\le1$.
  This implies that $\cT_n$ has at most $n$ ends.
\end{proof}

We endow $\cP$ with the weak topology:
$g_k\to g$ iff $g_k(\nu) \to g(\nu)$ for any $\nu\in\cVqm$. 
\begin{Lemma}\label{cvx}
  $\cP$ is a closed convex cone in
  $\cV_\mathrm{qm}^{\mathbf{R}}$ and
  $\cP_1=\{g\in\cP\ |\ g(\nu_\fm)=1\}$ is compact.
\end{Lemma}
\begin{proof}
  It is clear that if $g$ and $h$ are tree potentials, then so 
  is $ag+bh$ for any $a,b>0$. Let us show that $\cP$ is weakly
  closed.
  Condition~(P1) is easily seen to be preserved under pointwise
  limits. The same is true for~(P2)-(P3) by
  the following elementary fact: 
  pick a sequence $f_j$ of concave functions converging 
  pointwise towards $f$ on an open real interval. 
  Denote by $f'(x+)$ ($f'(x-)$) the right (left) derivative of 
  $f$ at $x$. 
  Then $f'(x+)\le\liminf f'_j(x+)\le\limsup f'_j(x-)\le f'(x-)$.
  
  As $\cV$ is compact, the set $\cF$ of functions on $\cV$ with values
  in $\Rbar$ is compact. 
  By~(P1), any $g\in\cP$ extends uniquely to an element of $\cF$.
  The above argument shows that the closure of $\cP$ in $\cF$ 
  is the union of $\cP$ and the function identically $+\infty$. 
  Hence $\cP_1$ is closed in $\cF$ and therefore weakly compact.
\end{proof}
\begin{Lemma}\label{L702}
  Any tree potential on $\cV$ is lower semicontinuous and restricts
  to a continuous tree potential on any finite subtree.
\end{Lemma}
\begin{proof}
  Condition~(P1) implies that the restriction of $g$ to any
  finite subtree $\cT\subset\cV$ is continuous.
  Denote by $g_\cT$ the function defined by $g(\nu)=g(\nu_0)$ where
  $\nu_0=\max[\nu_\fm,\nu]\cap\cT$. This is a tree potential on $\cV$
  which coincides with $g$ on $\cT$ and is continuous on $\cV$.
  By Lemma~\ref{lem-support} the support of
  $g$ is contained in the closure of an increasing union of finite
  subtrees $\cT_n$. Set $g_n=g_{\cT_n}$. 
  Then $g_n$ is continuous and increases pointwise to $g$ as $n\to\infty$. 
  Thus $g$ is lower semicontinuous.
\end{proof}                                
However, tree potentials are not necessarily continuous:
see Example~\ref{E701}.
\begin{Lemma}\label{L704}
  If $(g_i)_{i\in I}$ is any family of tree potentials, 
  then $g=\inf_ig_i$ is also a tree potential.
\end{Lemma}
\begin{proof}
  This is proved in the same way as the fact that the infimum of
  a family of concave functions on $\R$ is concave. The details are
  left to the reader.
\end{proof}
%
%%%%%%%%%%%%%%%%%%%%%%%%%%%%%%%%%%%%%%%%%%%%%%%%%%%%%%%%%%%%%%%%%%
%
\subsection{The tree Laplacian.}
Let us show how Borel measures give rise to tree potentials.
First pick any valuation $\mu\in\cV$ and
define $g_\mu:\cVqm\to[1,\infty)$ by 
$g_\mu(\nu)=\nu\cdot\mu$, where 
$\nu\cdot\mu=\a(\nu\wedge\mu)$ is the intersection product.
\begin{Lemma}\label{fundamental}
  The function $g_\mu$ is a tree potential.
\end{Lemma}
\begin{proof}
  Clearly $g_\mu$ is nonnegative, increasing and concave along
  totally ordered segments, and $g_\mu(\nu_\fm)=1$. 
  If $\vv\in T\nu_\fm$, then $D_\vv g_\mu=1$ if $\vv$ is 
  represented by $\mu$ and zero otherwise. 
  If $\vv\in T\nu$, $\nu\ne\mu$, 
  then $D_\vv g_\mu=-1$ if $\nu\ge\mu$ and 
  $\vv$ is represented by $\nu_\fm$, 
  one if $\vv$ is represented by $\mu$, and zero otherwise.
  This easily implies~(P2)-(P3).
\end{proof}
For an arbitrary positive measure $\rho\in\cM$ we set
\begin{equation}\label{eq-repr}
  g_\rho(\nu)
  =\int_\cV g_\mu(\nu)\,d\rho(\mu)
  =\int_\cV\mu\cdot\nu\,d\rho(\mu).
\end{equation}
As $\cV$ is compact, any Borel measure may be weakly approximated by 
a finite atomic measure (see \eg~\cite{Bou1}), 
hence $g_\rho$ is a tree potential in view of Lemma~\ref{cvx}.
We now describe the general identification of tree potentials with
positive measures on $\cV$. This identification is analogous
to the identification of  subharmonic functions (modulo
harmonic functions) with positive measures on $\R^n$.
\begin{Thm}{\cite[Theorems~7.61 and~7.64]{treeval}}\label{riesz}
  The map
  \begin{equation*}
    \cM\ni\rho 
    \mapsto g_\rho\in\cP
    \qquad
  \end{equation*}
  is a homeomorphism in the weak topology.
\end{Thm}
We refer to~\cite{treeval} for a proof of this result.
\begin{Def}
  The \emph{Laplacian} $\Delta g$ of a tree potential $g\in\cP$
  is by definition the unique positive Borel measure 
  such that $g_{\Delta g} = g$.
\end{Def}
\begin{Example}\label{E101}
  If $\phi\in R$ is irreducible, then 
  $\nu(\phi)=m(\phi)\nu\cdot\nu_\phi$, 
  hence $\nu\mapsto\nu(\phi)$ defines a tree potential
  whose Laplacian is $m(\phi)\nu_\phi$.
\end{Example}
The following properties characterize the
measure $\Delta g$: see Remark~\ref{R701}.
\begin{Prop}\label{prop-fdnt}
  Let $g\in \cP$. The measure $\rho=\Delta g$ is the unique positive
  measure on $\cV$ such that $\rho\,U(\vv)=D_\vv g$ for any tangent
  vector $\vv$ not represented by $\nu_\fm$; and
  $\rho\,U(\vv)=g(\nu_\fm)+D_\vv g$ when $\nu_\fm$ represents
  $\vv$. In particular, the total mass of $\Delta g$ is given by
  $g(\nu_\fm)$. Moreover, 
  \begin{equation*}
    \rho\{\nu_\fm\}=g_\rho(\nu_\fm)-\sum_{\vv \in T\nu_\fm}D_\vv g_\rho
    \qand
    \rho\{\nu\}=-\sum_{\vv\in T\nu}D_\vv g_\rho\ \text{for $\nu\ne\nu_\fm$}.
  \end{equation*}
\end{Prop}
The proof essentially follows by linearity from the case $\rho = \mu$,
and is left to the reader.  We shall also need:
\begin{Cor}\label{C901}
  If $g\in\cP$ and $\nu\in\cVqm$, then 
  $g(\nu)\le g(\nu_\fm)\a(\nu)$,
  with equality iff $\Delta g$ is supported 
  on $\{\mu\ge\nu\}$.
\end{Cor}
\begin{proof}
  We may assume $g(\nu_\fm)=\mass\,\Delta g=1$. 
  Then 
  $g(\nu)=1+\int_{\nu_\fm}^\nu\Delta g\{\mu'\ge\mu\}\,d\a(\mu)\le1+\a(\nu)-1$
  and equality holds iff $\Delta g$ is supported on $\{\mu\ge\nu\}$.
\end{proof}

\begin{Remark}\label{rem-extreme}
  The representation $g=\int g_\mu\,d\rho_g(\mu)$ 
  gives the Choquet decomposition
  of the tree potential $g$ in the closed convex cone $\cP$:
  the extremal points of $\cM$ are Dirac masses so 
  the extremal points of $\cP$ are of the form $g_\mu$.
\end{Remark}
\begin{Prop}
  Let $\cP'$ be the smallest closed positive subcone of 
  $\cV_\mathrm{qm}^{\mathbf{R}}$
  which is closed under infima and contains all
  functions of the form $\nu\mapsto\nu(\phi)$ 
  for $\phi\in R$ irreducible. Then $\cP'=\cP$.
\end{Prop}
\begin{proof}
  By Proposition~\ref{cvx}, Lemma~\ref{L704} and Example~\ref{E101} we
  have $\cP'\subset\cP$.  For the reverse inclusion it suffices to
  show that $g_\mu\in\cP'$ for $\mu\in\cV$.  As the set of divisorial
  valuations is dense in $\cV$, we may assume $\mu$ is divisorial.  But
  then we can find $\phi_1,\phi_2\in R$ irreducible such that
  $\mu=\nu_{\phi_1}\wedge\nu_{\phi_2}$.  Hence
  $g_\mu(\nu)=\min\{\nu(\phi_1),\nu(\phi_2)\}$, so $g_\mu\in\cP'$.
\end{proof}
\begin{Example}\label{E701}
  Let $\rho=\sum_{n\ge1}n^{-2}\nu_{y+nx}$.
  If $\nu_n=\nu_{y+nx,n^3}$, then $\nu_n\to\nu_\fm$ but
  $g_\rho(\nu_n)>n\to\infty>g_\rho(\nu_\fm)$. 
  This shows that tree potentials
  are not weakly continuous in general (see Lemma~\ref{L702}).
\end{Example}
%
%%%%%%%%%%%%%%%%%%%%%%%%%%%%%%%%%%%%%%%%%%%%%%%%%%%%%%%%%%%%%%%%%%
%
\subsection{Intersection of measures}
Using bilinearity we extend the intersection product on valuations
in $\cV$ to measures in $\cM$. More precisely, if $\rho,\sigma\in\cM$
then we define
\begin{equation*}
  \rho\cdot\sigma
  :=\iint_{\cV\times\cV}\mu\cdot\nu\,d\rho(\mu)d\sigma(\nu)
  =\int_{\cV}g_{\rho}(\nu)d\sigma(\nu)
  =\int_{\cV}g_{\sigma}(\mu)d\rho(\mu),
\end{equation*}
where $g_\sigma$ and $g_\rho$ are the tree potentials of $\rho$ and
$\sigma$, respectively. 
The last two equalities follow from Fubini and the definition
of $g_\rho$ and $g_\sigma$.

Since $\mu\cdot\nu\ge1$ for all $\mu,\nu$ we get
\begin{equation}\label{e910}
  \rho\cdot\sigma
  \ge\mass\rho\cdot\mass\sigma
  =g_\rho(\nu_\fm)g_\sigma(\nu_\fm)
\end{equation}
with equality iff the supports of $g_\rho$ and $g_\sigma$ intersect
only at $\nu_\fm$.
\begin{Prop}\label{P901}
  The intersection product is lower semicontinuous on $\cM$.
\end{Prop}
We refer to~\cite[Proposition~7.76]{treeval} 
for a proof of this result.  Note
that the intersection product is not continuous on $\cM$, and not even
on $\cV$ as exemplified by $\nu_n=\nu_{y-nx,2}$: here
$\nu_n\cdot\nu_n=2$ but $\nu_n\to\nu_\fm$ and $\nu_\fm\cdot\nu_\fm=1$.
%
%
%%%%%%%%%%%%%%%%%%%%%%%%%%%%%%%%%%%%%%%%%%%%%%%%%%%%%%%%%%%%%%%%%%%%%%%%%%%
%
%
\section{Tree transforms of psh functions}\label{sec-transform}
Fix a psh function $u$ and consider the real-valued
function $g_u$ on $\cVqm$ given by $g_u(\nu)=\nu(u)$, 
where $\nu(u)$ is defined by Theorem~\ref{thm-eval}.
We call $g_u$ the \emph{tree transform} of $u$.
We will show that $g_u$ is a tree potential. 
Its Laplacian $\rho_u=\Delta g_u$ is a measure on $\cV$
called the \emph{tree measure} of $u$ and contains a lot
of information on $u$.
We will try to understand what measures on $\cV$ arise in this way.
%
%%%%%%%%%%%%%%%%%%%%%%%%%%%%%%%%%%%%%%%%%%%%%%%%%%%%%%%%%%%%%%%%%%%%%%%%%%%
%
\subsection{Tree transforms are tree potentials}
Our main goal is to prove
\begin{Thm}\label{chartreetrans}
  The tree transform $g_u$ of any psh function $u$ is a tree potential
  on $\cV$ and the tree measure $\rho_u=\Delta g_u$
  has mass equal to the Lelong number of $u$.

  Moreover, $\rho_u$ puts no mass on formal 
  (\ie non-analytic)
  curve valuations and its mass on an analytic
  curve valuation $\nu_D$ is related to the mass of
  $dd^cu$ on $D$ as follows: $\rho_u\{\nu_D\}\ge \lambda\,m(D)$ 
  iff $dd^cu\ge \lambda[D]$.
\end{Thm}
Let us again emphasize that $\rho_u$ gives a very fine measurement
of the singularity of $u$ at $0$. We shall see in
Section~\ref{sec-inter}, that $\rho_u$ determines essentially all
generalized Lelong numbers of $u$ in the sense of Demailly. Here we
prove
\begin{Prop}\label{equalmeas}
  For two psh functions $u$ and $v$, the following three assertions
  are equivalent.
  \begin{enumerate}
  \item for all modifications $\pi:X\to(\C^2,0)$, and all
    points $p\in\pi^{-1}\{0\}$, we have
    $\nuL(\pi^*u,p)=\nuL(\pi^*v,p)$;
  \item $u$ and $v$ have the same tree transform: $g_u=g_v$;
  \item $u$ and $v$ have the same tree measure: $\rho_u=\rho_v$.
\end{enumerate}
\end{Prop}
\begin{Remark}
  It follows from Theorem~\ref{chartreetrans} that
  if $(x,y)$ are local coordinates, 
  then the function $t\mapsto\nu_{y,t}(u)$ is concave
  for $t\ge1$. Thus we recover the fact~\cite{kis2} that the
  Kiselman number $\nu^{x,y}_{a,b}$ is a concave function of 
  $(a,b)$.
\end{Remark}
\begin{Remark}\label{R902}
  If $u$ is psh with tree measure $\rho_u\in\cM$, then for any
  $\nu\in\cVqm$ quasimonomial we have
  $\nu(u)=\int_{\cV}\nu\cdot\mu\,d\rho_u(\mu)$.
  This follows from~\eqref{eq-repr}.
\end{Remark}
\begin{Remark}\label{R901}
  The tree transform $u\mapsto g_u$ inherits the main
  properties stated in Theorem~\ref{thm-eval}:
  compatibility, monotonicity, homogeneity,
  tropicality and semicontinuity.
%  If $u$ and $v$ are psh then $g_{u+v}=g_u+g_v$ and
%  $g_{\max\{u,v\}}=\min\{g_u,g_v\}$.
  In addition, if $u_n$ is the Demailly approximating 
  sequence of $u$, it follows from
  Proposition~\ref{P5} that $g_{u_n}\to g_u$ in $\cP$, hence
  $\rho_{u_n}\to\rho_u$ in $\cM$.
  These properties completely characterize the tree transform.
\end{Remark}
\begin{Example}\label{E401}
  The tree measure of $u=\log\max\{|x|,|y|\}$ is $\rho_u=\nu_\fm$,
  \ie a Dirac mass at the multiplicity valuation $\nu_\fm$.
  More generally, if $\phi$ is an irreducible germ of
  multiplicity $m=m(\phi)$, $x$ is a coordinate transverse
  to $x$ as in Section~\ref{def-qm}, and $1\le t\le\infty$,
  then the tree measure of $u=\log\max\{|\phi|^{1/m},|x|^t\}$
  is a Dirac mass at the valuation $\nu_{\phi,t}$.
  This follows from compatibility, homogeneity and
  tropicality.
\end{Example}
The proof of Theorem~\ref{chartreetrans}, given below, goes by
reduction to the algebraic case.  We define the \emph{tree transform}
$g_I$ of an ideal $I\subset R$ by $g_I(\nu)=\nu(I):=\min_{\phi\in
I}\nu(\phi)$.
\begin{Prop}\label{Plog}
  The tree transform $g_I$ of any ideal $I\subset R$ is a tree potential.
  Its Laplacian $\rho_I=\Delta g_I$ has mass $m(I)\= \nu_\fm (I)$ and is
  an atomic measure supported on finitely many 
  divisorial and (analytic) curve valuations.
\end{Prop}
\begin{Remark}
  In~\cite[Theorem 8.2]{treeval} we characterize measures on $\cV$ of
  the form $\rho_I$: they are atomic measures whose mass is a multiple
  of the generic multiplicity at any divisorial valuation, and a
  multiple of the multiplicity at any curve valuation.
\end{Remark}
\begin{proof}
  For $\phi\in\fm$, set $g_\phi(\nu)=\nu(\phi)$. 
  When $\phi$ is irreducible, Example~\ref{E101}
  shows that $g_\phi$ is a tree potential 
  and that $g_\phi$ is a piecewise affine function with
  integer slopes on any segment in $\cV$ parameterized by skewness.
  By additivity in $\cP$ and unique
  factorization in $R$, the same properties hold when $\phi$ is reducible. 
  
  Let $S\subset I$ be a finite set of generators for $I$.  Then
  $g_I=\min_{\phi\in S}g_\phi$. Thanks to Lemma~\ref{L704}, $g_I$ is a
  tree potential. By Proposition~\ref{prop-fdnt}, the mass of 
  $\rho_I:=\Delta g_I$ is given by $g_I(\nu_\fm)=m(I)$.

  It is clear that $g_I$ is supported on the smallest subtree of $\cV$
  containing $\nu_\fm$ and any $\nu_\psi$, where $\psi$ ranges over
  the irreducible factors of the elements of $S$.  This is a finite
  subtree $\cS$. Moreover, it follows from the preceding computation
  that on any segment in $\cS$ parameterized by skewness, $g_I$ is a
  piecewise affine function with integer slopes. Thus $\rho_I$ is a
  finite sum of point masses, taken over valuations that are either
  ends or branch points in $\cS$, or regular points in $\cS$
  where $g_I$ fails to be locally affine. From the integer slope
  property we conclude that $\rho_I=\sum_{i=1}^rn_i\nu_i$, where
  $\nu_i$ are divisorial (\ie have rational skewness) or curve
  valuations and $n_i$ are positive integers.
\end{proof}
\begin{proof}[Proof of Theorem~\ref{chartreetrans}]
  First suppose $u$ has logarithmic singularities, and write
  $u=\frac{c}2\log\sum_{i=1}^n|\phi_i|^2$ for holomorphic $\phi_i$ and
  $c>0$.  Then $g_u(\nu)=c\,g_I(\nu)= c \min \nu(\phi_i)$, where $I$
  is the ideal generated by the $\phi_i$ (see
  Proposition~\ref{P4}). Hence $g_u$ is a tree potential in this case
  by Proposition~\ref{Plog}. Further, the mass of $\Delta g_u$ equals
  $m(I)$, which is the Lelong number of $u$. 

  In the general case, we use Proposition~\ref{P5}.  Let $u_n$ be a
  sequence of psh functions with logarithmic singularities such that
  $u_n\to u$, and $g_{u_n}\to g_u$ pointwise.  Since $\cP$ is weakly
  closed (Lemma~\ref{cvx}) it follows that $g_u$ is a tree
  potential. That $\Delta g_u$ has the right mass follows since
  $g_{u_n}(\nu_\fm)\to g_u(\nu_\fm)$.

  \smallskip
  For the second assertion, 
  first consider an irreducible holomorphic germ $\phi\in\fm$,
  and suppose $dd^cu\ge\lambda[D]$, where $D=\{\phi=0\}$ and $\lambda>0$.
  Then $u\le\lambda\log|\phi|$ and 
  $g_u(\nu_{\phi,t})\ge\lambda\,m(\phi)t$ for all $t\ge1$. 
  This implies $\rho_u\{\nu_\phi\}\ge\lambda\,m(\phi)$. 

  Conversely, suppose $\phi\in\hat{\fm}$ is a formal germ
  and $\rho_u\{\nu_\phi\}\ge\lambda\,m(\phi)$ with $\lambda>0$.
  Fix any $\e>0$.
  Then $g_u(\nu_{\phi,t})\ge(\lambda-\e)\,m(\phi)t$ for large $t$.

  First assume $u=\frac{c}2\log\sum_{i=1}^n|\phi_i|^2$ has
  logarithmic singularities.
  Then $g_{\phi_i}(\nu_{\phi,t})\ge c^{-1}(\lambda-\e)\,m(\phi)t$ for all $i$,
  which implies that $\phi^k$ divides $\phi_i$ for all $i$, where
  $k$ is the smallest integer larger than $c^{-1}\lambda$.
  This implies that the curve $D=\{\phi=0\}$ is analytic,
  and that $u\le ck\log|\phi|+O(1)$.
  Hence $dd^cu\ge ck[D]\ge\lambda[D]$.

  In the general case we have
  $|g_{u_n}(\nu_{\phi,t})-g_u(\nu_{\phi,t})|
  \le A(\nu_{\phi,t})/n\le m(\phi)t/n$ 
  by~\eqref{e904} and the estimates from Section~\ref{thinness}.
  Hence, for $n$ large, 
  $g_{u_n}(\nu_{\phi,t})\ge(\lambda-2\e)\,m(\phi)t$ 
  for large $t$.
  By what precedes $D=\{\phi=0\}$ is analytic, 
  $dd^cu_n$ puts mass at least $\lambda$ on $D$,
  and $u_n\le\lambda\log|\phi|+O(1)$.
  By~\eqref{e802}, $u\le\lambda\log|\phi|+O(1)$, 
  so $dd^cu\ge\lambda[D]$.
\end{proof}
\begin{proof}[Proof of Proposition~\ref{equalmeas}]
  The equivalence of $(2)$ and $(3)$ is a consequence of
  Theorem~\ref{riesz}. The implication $(2)\Rightarrow (1)$ follows
  from Proposition~\ref{P-interpret}, and the following fact. The
  Lelong number at a point $p$ equals the Lelong number at a generic
  point on the exceptional divisor obtained by blowing up $p$. 
  Finally, suppose $(1)$ is true. Then
  Proposition~\ref{P-interpret} shows that $g_u(\nu) = g_v(\nu)$ for
  all divisorial valuations $\nu$. As tree potentials are continuous on
  finite subtrees by Lemma~\ref{L702}, and divisorial valuations are
  dense on any finite subtree in $\cVqm$, we conclude $g_u = g_v$.
\end{proof}

%
%%%%%%%%%%%%%%%%%%%%%%%%%%%%%%%%%%%%%%%%%%%%%%%%%%%%%%%%%%%%%%%%%%%%%%%%%%%
%
\subsection{Representation of measures by psh functions}
% When $\rho=\rho_u\in\cM$ is the tree measure of a psh function $u$, 
% we say that $\rho$ is \emph{represented} by $u$. 
% Such a measure cannot put mass on any formal 
% curve valuations by Theorem~\ref{chartreetrans}.
% In general it seems quite hard to characterize 
% the measures in $\cM$ that are represented by psh
% functions.
% Less ambitiously we may ask when a Dirac mass 
% at a valuation $\nu$ is represented by a psh function. 
% A curve valuation is represented by a psh function iff 
% the curve is analytic. If $\nu$ is quasimonomial of multiplicity $m$, 
% say $\nu=\nu_{\phi,t}$, where $m(\phi)=m$, then
% Remark~\ref{R901} shows that $\nu$ is represented by 
% $u=\log\max\{|x|^{t},|\phi|^{1/m}\}$, for any choice 
% of transverse curve $\{x=0\}$.

When $\rho=\rho_u\in\cM$ is the tree measure of a 
psh function $u$, 
we say that $\rho$ is \emph{represented} by $u$. 
In general it seems quite hard to characterize 
the measures in $\cM$ that are represented by psh
functions (or currents). 
On the one hand, Example~\ref{E401} shows that
any Dirac mass at a quasimonomial or analytic curve valuation
is represented by a psh function. 
By taking sums and limits we obtain many more measures:
see Example~\ref{E402} for an interesting example.

On the other hand, there are also some restrictions.
The tree measure of a psh function cannot put mass on 
any formal curve valuation by Theorem~\ref{chartreetrans}.
There are likely some restrictions at infinitely singular 
valuations, too, but we do have
\begin{Prop}\label{repr-solenoid}
  Let $\nu$ be an infinitely singular valuation, given by a Puiseux
  series $\hat{\phi}=\sum_1^\infty a_jx^{\hb_j}$ as in
  Section~\ref{puiseux}.  Suppose there exists $r>0$ such that
  $\sum_1^\infty|a_j|^2\hb_jr^{2\hb_j}<\infty$.  Then $\nu$ is
  represented by a psh function.
\end{Prop}
With stronger assumptions, 
one should be able to analyze more precisely the set of
psh functions (or currents) representing $\nu$. One may for instance
ask for conditions on the $a_j$'s ensuring the existence of a
\emph{unique} extremal positive closed $(1,1)$-current representing $\nu$.
We refer to~\cite{slod} for related problems and
to~\cite[Section~6]{kiwi} for results in this direction.
\begin{proof}
  We may assume that $(x,y)$ in the definition of the Puiseux series
  are global coordinates on $\C^2$.
  Let $\hat{\phi}_n=\sum_1^na_jx^{\hb_j}$ as above, and let
  $\phi_n\in R$ be the minimal polynomial of $\hat{\phi}_n$.
  Write $m_n=m(\phi_n)$. The divisorial valuation 
  $\nu_k=\nu_{\phi_n}\wedge\nu_{\phi_k}$ is independent of $n$
  for $n>k$ and $\nu_k$ increases to $\nu$ as $k\to\infty$.

  Pick $r_0>0$ such that $\sum_1^\infty|a_j|^2\hb_jr^{2\hb_j}<\infty$,
  and define $C=\sum_1^\infty|a_j|^2\hb_jr_0^{2(\hb_j-\hb_1)}$.
  Let $B_r=\Delta_r\times\Delta_1$ for $r\in(0,r_0]$.
  Define $u_n=\frac{1}{m_n}\log|\phi_n|$ and $T_n=dd^cu_n$.
  Then $T_n$ forms a sequence of positive closed currents on $B_{r_0}$.
  Further, $T_n$ is the pushforward of the current of integration
  on the disk $\Delta_n=\{|t|<r_0^{1/m_n}\}$ by the map
  $\psi_n(t)=(t^{m_n},\sum_1^na_jt^{m_n\hb_j})$.
  Notice that $\psi_n(\Delta_n)\in\{|x|<r_0, |y|\le C|x|^{\hb_1}\}$.

  Let $b_n(r)$ be the mass of (the trace measure of) $T_n$ in $B_r$.
  We have $b_n(r)=\pi(r^2+\sum|a_j|^2\hb_j r^{2\hb_j})$, 
  hence $b_n(r)\le \pi(r^2+Cr^{2\hb_1})$ for $r<r_0$.  
  Thus we may extract a subsequence $T_{n_j}$ that
  converges to a positive closed current $T$ in $B_{r_0}$. The mass
  $b(r)$ of $T$ in $B_r$ also satisfies $b(r)\le \pi(r^2+Cr^{2\hb_1})$.  As
  $\hb_1>1$, this implies that the Lelong number of $T$ is at most 1
  (see~\ref{eq-def-lelong}).

  On the other hand, for $n\ge k$ we have
  $\nu_k(T_n)=\nu_k(\phi_n)/m_n=\alpha(\nu_k)$.  By semicontinuity
  this gives $\nu_k(T)\ge\alpha(\nu_k)\ge\alpha(\nu_k)\nu_\fm(T)$.
  Let $\rho=\rho_T$ be the measure represented by $T$.  By
  Corollary~\ref{C901}, $\rho$ has mass 1 and is supported on
  $\{\mu\ge\nu_k\}$ for all $k$.  Thus $\rho$ is a point mass at
  $\nu$, completing the proof.
\end{proof}
%
%%%%%%%%%%%%%%%%%%%%%%%%%%%%%%%%%%%%%%%%%%%%%%%%%%%%%%%%%%%%%%%%%%%%%%%%%%%
%
%
\section{Attenuation of singularities of currents}\label{sec-att}
\begin{figure}[ht]
  \includegraphics[width=0.6\textwidth]{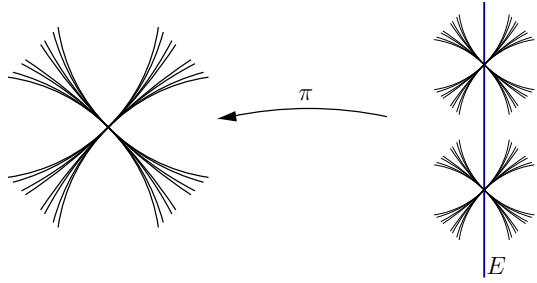}
  \caption{The Cantor bouquet of curves splits into two
    isomorphic parts after a single blowup of the origin.
    See Example~\ref{E402}.}\label{F5}
\end{figure}

Positive closed $(1,1)$-currents in many ways generalize curves
(in dimension 2). Here we shall prove a theorem that 
generalizes embedded resolution of plane curve singularities.
\begin{Thm}\label{thm-attlocal}
  Let $T$ be any positive closed $(1,1)$-current near the origin
  in $\C^2$. Then for any $\eta,\e>0$ there exists a 
  modification $\pi:X\to(\C^2,0)$ 
  and positive closed   
  $(1,1)$-currents $S_1$, $S_2$ on $X$,
  such that $\pi^*T=S_1+S_2$ and:
  \begin{itemize}
  \item
    the support of $S_1$ is a curve with normal crossing singularities;
  \item 
    $\sup_{p\in\pi^{-1}(0)}\nuL(S_2,p)^{1+\eta} 
    \le\sum_{p\in\pi^{-1}(0)}\nuL(S_2,p)^{1+\eta}\le\e$.
  \end{itemize}
\end{Thm}
This statement implies the corresponding global statement, when $T$ is
defined on a compact complex surface. 
Thus Theorem~\ref{thm-attlocal} strengthens the main result 
of~\cite{guedj} in two ways, as the 
method there---borrowed from~\cite{Mim}---only gives the 
weaker bound $\sum_{\pi^{-1}(0)}\nuL(S_2,p)^2\le\e$. 
We also recover the following result by
Mimouni~\cite[Th{\'e}or{\`e}me~III.1.2]{Mim}.  See below for a definition of
strict transform.
\begin{Cor}\label{cor-mim}
  Suppose $T$ does not charge any curve, and fix $\e>0$. 
  Then there exists a modification $\pi:X\to (\C^2,0)$ 
  such that the strict transform of $T$ by $\pi$ has all its 
  Lelong numbers bounded by $\e$.
\end{Cor}
\begin{Remark}\label{R711}
  Theorem~\ref{thm-attlocal} fails for $\eta=0$ in general. 
  As the proof of Lemma~\ref{lemma-key} below shows, 
  any psh function $T$ whose tree measure $\rho_T$ 
  has no atoms and is supported on the set of smooth,
  analytic, curve valuations, yields a counterexample.
  See Remark~\ref{finalrem} for more details and Example~\ref{E402}
  for an explicit construction.
\end{Remark}
\begin{Example}\label{E402}
  Let $\Sigma\=\{-1,+1\}^{\N^*}$, and $\rho$ be the 
  uniformly distributed measure on $\Sigma$, \ie the product measure
  $\rho=\otimes_1^\infty\rho_i$, where $\rho_i\{+1\}=\rho_i\{-1\}=1/2$
  for any $i$.
  For $\sigma=(\sigma_i)_1^\infty\in\Sigma$ set
  $f_\sigma(x)\=\sum_{i\ge1}\sigma_i x^i$.
  Now define
  \begin{equation*}
    T=\int_{\Sigma}[y=f_\sigma(x)]\,d\rho(\sigma),
  \end{equation*}
  where $[y=f_\sigma(x)]$ denotes the current of integration 
  on the curve $\{y=f_\sigma(x)\}$. 
  The set $\Sigma$ can be thought of
  as a Cantor set on two symbols, and the collection 
  $\{y=f_\sigma(x)\}_{\sigma\in\Sigma}$,
  as a ``Cantor bouquet'' of smooth curves.
  Indeed, by sending $\sigma\in\Sigma$ to the 
  curve valuation $\nu_{y-f_\sigma(x)}$ 
  we obtain a homeomorphism of $\Sigma$ (with the product topology)
  onto a compact subset of $\cV$.
  The tree measure 
  $\rho_T$ of $T$ on $\cV$ is then the pushforward 
  of $\rho$ on $\Sigma$.

  Under any modification $\pi$, the strict transform of the bouquet
  will be a finite union of bouquets all of which are isomorphic to
  the original one. See Figure~\ref{F5}.
  The sum of the Lelong numbers of the strict
  transform of $T$ will always equal one, the Lelong number of $T$
  at the origin.
\end{Example}
\begin{Remark}\label{R712}
  Let us note a strong similarity between the statement of
  Theorem~\ref{thm-attlocal}
  and the definition of Zariski's decomposition as in
  ~\cite[Theorem~3.1]{BDP} (see also~\cite{Bou}). 
%  It would be
%  interesting to clarify the interplay between our local statement and
%  its global analog.
\end{Remark}
\begin{proof}[Proof of Theorem~\ref{thm-attlocal}] 
%   We decompose the proof into two steps. First we prove that we may
%   pick $\pi$ such that $\pi^*T = S_1 + S_2$ as in the statement of the
%   theorem with $\sup \nuL (S_2,p)\le\e$.  Then we show the more
%   precise bound on the sum of powers of Lelong numbers of $S_2$.

  We shall prove that we may pick $\pi$ such that 
  $\pi^*T=S_1+S_2$ as in the statement of the theorem with 
  $\sup_p\nuL(S_2,p)$ arbitrarily small and $\sum_p\nuL(S_2,p)$ 
  uniformly bounded. This will make the
  $\ell^{1+\eta}$-norm arbitrarily small.
  
  By Siu's Theorem, $T=T_1+T_2$, where $T_1=\sum_{j=1}^na_j[D_j]$ with
  $a_j>0$ and $D_j$ irreducible, and $T_2\not\ge\e[D]$ for all
  irreducible curves $D$.

  Denote by $\Gast_\pi$ the set of exceptional components of 
  a modification $\pi$.
  The pull-back $\pi^*T$ is a positive closed $(1,1)$-current which
  charges any curve $E\in\Gast_\pi$; more precisely
  $\pi^*T\ge\div_E(\pi^*T)[E]$ for any $E\in\Gast_\pi$. By
  Proposition~\ref{P-interpret}, $\div_E(\pi^*T)$ is also the Lelong
  number of $\pi^*T$ at a generic point $p\in E$. The \emph{strict
    transform} of $T$ by $\pi$ is by definition
  \begin{equation*}
    T'\=\pi^*T -\sum_{E\in\Gast_\pi}\div_E(\pi^*T)[E].
  \end{equation*}
  It is a positive closed $(1,1)$-current with zero Lelong number except
  on a countable subset of $\pi^{-1}(0)$.
  
  We now prove that we may pick $\pi$ such that the Lelong numbers of
  the strict transform of $T_2$ are arbitrarily small. As we shall
  see, composing $\pi$ with further blowups does not destroy the
  latter property.  
  Hence we may assume that the total transform of
  the curve $\bigcup D_j$ has simple normal crossings, so that the theorem
  holds with $S_2$ being the strict transform of $T_2$ and
  $S_1=\pi^*T-S_2$.

  By this argument, we are reduced to the case $T_1=0$.
  In other words, we may assume that 
  $T\not\ge\e[D]$ for all irreducible curves $D$.
  We shall pick $S_2=T'$, the strict transform of $T$.

  Let $\rho_T$ be the tree measure of $T$. The assumption
  on $T$ means that $\rho_T\{\nu_D\}<\e\,m(D)$
  for all $D$: see Theorem~\ref{chartreetrans}.
  Our aim is to estimate the Lelong numbers of $T'$
  in terms of the mass of $\rho_T$ on particular 
  subregions in the valuative tree.
  To do so, we rely in an essential way on the 
  results described in Section~\ref{divisorial}.
  
  Consider any modification $\pi:X\to(\C^2,0)$ and pick a point
  $p\in\pi^{-1}(0)$.  Denote by $\nu_p$ the divisorial valuation
  associated to the blowup of $p$, and define an open subset
  $U(p)\subset\cV$ as follows (see Figure~\ref{F1} on p.\pageref{F1}
  and compare with Proposition~\ref{prop-divisorial}):
  \begin{itemize}
  \item
    If $p$ is a regular point on $\pi^{-1}(0)$, lying on a
    unique exceptional component $E$, then $\nu_p>\nu_E$ and 
    $U(p):=\{\mu \ ; \ \mu\wedge\nu_p> \nu_E \}$.
  \item
    If $p$ is a singular point on $\pi^{-1}(0)$, then $p=E\cap E'$,
    where $\nu_E<\nu_p<\nu_{E'}$, and
    $U(p):=\{\mu;\ \nu_E<\mu\wedge\nu_{E'}<\nu_{E'}\}$.
  \end{itemize}
  In can be shown that $U(p)$ is exactly the set of valuations in
  $\cV$ whose center on $X$ is the point $p$:
  see~\cite[Proposition~6.32]{treeval}.  
  For fixed $\pi$, the sets $U(p)$ form a
  disjoint open cover of $\cV\setminus\{\nu_E\ ;\  E\in\Gast_\pi\}$. 
  Recall that $\nuL(T',p)$ denotes the Lelong number of
  the current $T'$ at $p$.  We shall prove:
  \begin{Lemma}\label{lemma-key}
    For any current $T$, any modification 
    $\pi:X\to(\C^2,0)$ and any point $p\in\pi^{-1}(0)$ we have 
    $\rho_T(U(p))\ge b(\nu_p)\,\nuL(T',p)$.
  \end{Lemma}
  \begin{Lemma}\label{L102}
    Let $\rho\in\cM$ be a positive measure on $\cV$ such that
    $\rho\{\nu_D\}<\e\,m(D)$ for every irreducible
    curve $D$. Then there exists a modification $\pi:X\to(\C^2,0)$ 
    such that $\rho(U(p))\le\e\,b(\nu_p)$ for every 
    $p\in\pi^{-1}(0)$.
  \end{Lemma}
  In view of the reductions above, we obtain, for any $\e>0$,
  the existence of $\pi$ such that $\nuL(T',p)\le\e$ 
  for all $p\in\pi^{-1}(0)$. On the other hand, since 
  $b(\nu_p)\ge1$ for all $p$ and the sets $U(p)$ are
  disjoint, Lemma~\ref{lemma-key}
  shows that $\sum_{p\in\pi^{-1}(0)}\nuL(T',p)$ is
  uniformly bounded by the mass of $\rho_T$, \ie 
  the Lelong number $\nuL(T)$ of $T$ at the origin.
  For $\eta>0$, we get 
  $\sum_{p\in\pi^{-1}(0)}\nuL(T',p)^{1+\eta}\le\e^\eta\nuL(T)$.
  If $\eta$ and $T$ are fixed, we can then make $\e^\eta\nuL(T)$ 
  arbitrarily small.  This concludes
  the proof of Theorem~\ref{thm-attlocal}.
\end{proof}
\begin{proof}[Proof of Lemma~\ref{lemma-key}]
  To simplify notation, we shall
  write $b_E$ for $b(\nu_E)$, $\a_E$ for $\a(\nu_E)$ etc.
  Let $\mu_p$ be the multiplicity valuation at $p$ in $X$
  and write $U=U(p)$.
  
  First suppose $p\in E$ is a regular point of $\pi^{-1}(0)$.
  By definition $T'=\pi^*T-\div_E(\pi^*T)[E]$. 
  In particular $\nuL(T',p)=\mu_p(\pi^*T)-\div_E(\pi^*T)$. 
  We have $\nu_E=b_E^{-1}\pi_*\div_E$ and $\nu_p=b_p^{-1}\pi_*\mu_p$.
  Moreover, $b_p=b_E$ by Proposition~\ref{prop-divisorial}~(i), so
  \begin{equation*}
    \nuL(T',p)=b_p(\nu_p(T)-\nu_E(T)).
  \end{equation*}
  By Lemma~\ref{lem-estim} we have $\a_p-\a_E=b_p^{-2}$.
  Using Remark~\ref{R902} we have
  \begin{align*}
    \nu_p(T)-\nu_E(T) 
    &=\int_\cV\left(\nu_p\cdot\nu-\nu_E\cdot\nu\right)\,d\rho_T(\nu)
    =\int_U\left(\nu_p\cdot\nu-\nu_E\cdot\nu\right)\,d\rho_T(\nu)\\
    &\le(\a_p-\a_E)\,\rho_T(U)=b_p^{-2}\,\rho_T(U).
  \end{align*}  
  This concludes the proof in this case.

  \smallskip
  Now suppose $p=E\cap E'$ is a singular point of $\pi^{-1}(0)$.
  Then $T'=\pi^*T-\div_E(\pi^*T)[E]-\div_{E'}(\pi^*T)[E']$.  
  We have $\nu_E=b_E^{-1}\pi_*\div_E$, 
  $\nu_{E'}= b_{E'}^{-1}\pi_*\div_{E'}$
  and $\nu_p=b_p^{-1}\pi_*\nu_p$. 
  Moreover, $b_p=b_E+b_{E'}$ by 
  Proposition~\ref{prop-divisorial}~(ii), so
  \begin{equation*}
    \nuL(T',p)
    =(b_E+b_{E'})\,\nu_p(T)-b_E\nu_E(T)-b_{E'}\nu_{E'}(T).
  \end{equation*}
  We may assume that $\nu_E<\nu_p<\nu_{E'}$.
  By Lemma~\ref{lem-estim} we then have
  $\a_p-\a_E=b_E^{-1}(b_E+b_{E'})^{-1}$
  and $\a_{E'}-\a_p=b_{E'}^{-1}(b_E+b_{E'})^{-1}$.
  This implies that
  \begin{equation*}
    \nuL(T',p)
    =\frac1{b_p}\left(
      \frac{\nu_p(T)-\nu_E(T)}{\a_p-\a_E}-
      \frac{\nu_{E'}(T)-\nu_p(T)}{\a_{E'}-\a_p}
    \right).
  \end{equation*}
  As above we have 
  $\nu_p(T)-\nu_E(T)\le(\a_p-\a_E)\,\rho_T\{\mu\wedge\nu_p>\nu_E\}$.
  Furthermore,
  \begin{multline*}
    \nu_{E'}(T)-\nu_p(T)
    =\int_\cV\left(\nu_p\cdot\nu-\nu_E\cdot\nu\right)\,d\rho_T(\nu)\\
    \ge\int_{\nu\ge\nu_{E'}}
    \left(\nu_p\cdot\nu-\nu_E\cdot\nu\right)\,d\rho_T(\nu)
    \ge\rho_T\{\mu\ge\nu_{E'}\}(\a_{E'} - \a_p),
  \end{multline*}
  where the first equality follows from Remark~\ref{R902}.
  Thus we conclude
  \begin{equation*}
    \nuL(T',p) \le\frac{1}{b_p} 
    \left(
      \rho_T\{\mu\wedge\nu_{F}>\nu_E\}-\rho_T\{ \mu\ge\nu_{E'}\}
    \right) 
    =\frac{1}{b_p}\rho_T(U),
  \end{equation*}
  which completes the proof of Lemma~\ref{lemma-key}.
\end{proof} 
\begin{Remark}\label{finalrem}
  The proof shows that equality holds, \ie
  $\rho_T(U(p))=b(\nu_p)\nuL(T',p)$, when $\rho_T$ is supported on the
  set of smooth curve valuation.  Compare
  Remark~\ref{R711}.
\end{Remark}
\begin{proof}[Proof of Lemma~\ref{L102}]
  Define $\cT=\{\nu \ | \ \rho\{\mu\ge\nu\} \ge\e\,m(\nu)\}$.
  It is clear that $\nu\in\cT$, $\nu'\le\nu$ implies $\nu'\in\cT$, 
  hence $\cT$ is a subtree of $\cV$.
  It is moreover a finite subtree, with at most
  $\e^{-1}\times\mass\rho$ ends.
  Our assumption implies that all these ends 
  are quasimonomial valuations.

  Now pick a modification $\pi:X\to(\C^2,0)$
  with the following properties:
  \begin{enumerate}
  \item[(i)]
    every end in $\cT$ is dominated by $\nu_E$, for some
    exceptional component $E$;
  \item[(ii)]
    whenever $E$, $E'$ are exceptional components of $\pi$ that 
    intersect in $X$, we have $b(E)+b(E')>\mass\rho/\e$.
  \end{enumerate}
  We may achieve~(i) as each end in $\cT$ are 
  quasimonomial, hence dominated by some 
  divisorial valuation, which we may assume to be of the form
  $\nu_E$ for some exceptional component $E$ of $\pi$.
  If~(ii) would fail for some pair $E$, $E'$, 
  then we may compose $\pi$ with 
  the blowup at $E\cap E'$. This creates a new exceptional
  component $F$ with $b(F)>\max\{b(E),b(E')\}$.
  Thus~(ii) holds after finitely many further blowups.

  Now pick any $p\in\pi^{-1}(0)$. 
  If $p=E\cap E'$ is a singular point on $\pi^{-1}(0)$,
  then the conclusion of Lemma~\ref{L102} is immediate
  as $b(F)=b(E)+b(E')>\mass\rho/\e$.

  Hence suppose $p$ is a regular point on $\pi^{-1}(0)$,
  belonging to a unique exceptional component $E$.
  Then $\nu_p$ does not represent the same tangent vector
  as any $\nu_E'$, $E'$ ranging over exceptional components
  of $\pi$. By~(i) this implies that $U(p)\cap\cT=\emptyset$.
  For any $\nu\in\,]\nu_E,\nu_p]$ we thus have 
  $\rho\{\mu\ge\nu\}<\e\,m(\nu)=\e\,b(\nu_p)$.
  As $\nu\to\nu_E$ we conclude
  $\rho(U(p))\le\e\,b(\nu_p)$, which completes the proof.
\end{proof}
%
%
%%%%%%%%%%%%%%%%%%%%%%%%%%%%%%%%%%%%%%%%%%%%%%%%%%%%%%%%%%%%%%%%%%%%%%%%%%%
%
%
\section{Intersection formula}\label{sec-inter}
Our aim is to relate the mass at the origin of the intersection
product of two positive closed $(1,1)$ currents to their tree
measures. An optimistic guess is
\begin{equation}\label{eq-inter}
  dd^cu\wedge dd^c v\,\{0\}
  =\rho_u\cdot\rho_v
%  =\iint\limits_{\cV\times\cV}\mu\cdot\nu\,d\rho_u(\mu)d\rho_v(\nu),
  =\iint_{\cV\times\cV}\mu\cdot\nu\,d\rho_u(\mu)d\rho_v(\nu),
\end{equation}
where $\rho_u$, $\rho_v$ are the tree measures of 
$u$ and $v$, respectively.
However, this is indeed too much to hope for.
\begin{Example}\label{example}
  Let $u=\max\{-\sqrt{-\log|x|},\log|y|\}$, and $v=\log|x|$.  
  Then $dd^cu\wedge dd^cv\,\{0\}$ equals the mass of $dd^cu|_{x=0}$ which
  is one. On the other hand, the Lelong number of $u$ is zero,
  hence $\rho_u=0$ so $\rho_u\cdot\rho_v=0$ and~\eqref{eq-inter} fails.
\end{Example}
We may more precisely conjecture~\eqref{eq-inter} as soon as 
neither current charges an analytic curve.  Note that this would in
particular imply that the admissible wedge product of two currents
with \emph{zero} Lelong number never charges the origin---something
that seems quite hard to prove.

Our aim in this section is to give partial results in the
direction of~\eqref{eq-inter}. We prove that equality holds when 
either $u$ or $v$ has logarithmic singularities 
(Proposition~\ref{p-inter}). 
We prove that a lower bound always holds
(Proposition~\ref{p-lower}). Finally we show that equality
holds whenever $v$ is a psh weight for which $e^v$ is
\Holder continuous (Proposition~\ref{p-upper}). This
gives an interpretation of a generalized Lelong number as an average
of valuations.
%
%%%%%%%%%%%%%%%%%%%%%%%%%%%%%%%%%%%%%%%%%%%%%%%%%%%%%%%%%%%%%%%%%%%%%
% 
\subsection{Psh functions with logarithmic singularities}
Here we prove
\begin{Prop}\label{p-inter}
  Let $u$ and $v$ be psh functions such that $dd^cu\wedge dd^cv$ is
  admissible.
  Suppose we are in one of the following two cases:
  \begin{itemize}
  \item either  $u$ is an arbitrary psh function, $v$ has 
    logarithmic singularities and $dd^cv$ does not charge any curve;
  \item or $u,v$ both have logarithmic singularities.
  \end{itemize}
  Then~\eqref{eq-inter} holds.
\end{Prop}
\begin{Remark}
  An algebraic version of this result can be found 
  in~\cite[Section~8.1.4]{treeval}.
  Pick $I,J$ two analytic ideals, generated by finitely many
  $f_i$'s and $g_j$'s respectively and define 
  $u=\frac12\log\sum|f_i|^2$, 
  $v=\frac12\log\sum|g_j|^2$. 
  Then the mass $dd^cu\wedge dd^cv\{0\}$ 
  can be naturally interpreted as a mixed
  multiplicity of the ideals $I$ and $J$. 
  Theorem~8.13 from~\cite{treeval} asserts
  that in a quite general algebraic setting,
  the mixed multiplicity of two ideals can be computed using their tree
  transforms.
\end{Remark}
\begin{proof}[Proof of Proposition~\ref{p-inter}]
  First assume $v$ has logarithmic
  singularities and does not charge any curves.  The tree transform of
  $v$ is a sum of Dirac masses at finitely many divisorial valuations,
  $\rho_u=\sum_1^k c_i\nu_i$.
  \begin{Lemma}\label{lem-log-sing}
    Let $v$, $v'$ be two psh functions
    with logarithmic singularities, whose tree transforms coincide. Then
    the difference $v-v'$ is a bounded function.
  \end{Lemma} 
  We hence choose for each $1\le i\le k$, an irreducible $\phi_i\in\fm$,
  and $t_i\ge 1$ so that $\nu_i=\nu_{\phi_i,t_i}$.  By the preceding
  lemma, and Proposition~\ref{p-semi} we may replace $v$ by the sum
  $\sum_{i=1}^k c_i \log\max\{\|p\|^{t_i},\ |\phi_i|^{1/m_i}\}$,
  $m_i=m(\phi_i)$.  By linearity of both sides of the
  equation~\eqref{eq-inter}, we are reduced to the case where the
  measure $\rho_v$ is supported at a single valuation, i.e.  we can
  suppose $k=1$, and $ v=\log\max\{\|p\|^t, |\phi|^{1/m}\} $, for an
  irreducible $\phi\in\fm$, and $t\ge1$. Now $dd^cu\wedge
  dd^cv\,\{0\}=\nu_{\phi,t}(u)$ follows from Proposition~\ref{P6}.
  This concludes the proof in this case.

  When $u$ and $v$ both have logarithmic singularities, we write 
  $u=u'+\sum a_i\log|\phi_i|$, 
  $v=v'+\sum b_j\log|\psi_j|$, where 
  $dd^cu'$ and $dd^cv'$ do not charge any curve, 
  $\phi_i,\psi_j\in\fm$ are irreducible analytic functions and $a_i,b_j>0$.  
  By linearity and from
  what precedes we are reduced to the case $u=\log|\phi|$,
  $v=\log|\psi|$ for distinct irreducible $\phi,\psi\in\fm$.  In this
  case, $dd^cu\wedge dd^cv$ is the pull-back under the finite map
  $F(z,w)=(\phi(z,w),\psi(z,w))$ of the measure 
  $\mu=dd^c\log|z|\wedge dd^c\log|w|$. 
  The mass of $\mu$ at $0$ is one, and the topological
  degree of $F$ is exactly the intersection product of the 
  two curves $C=\{\phi=0\}$ and $D=\{\psi=0\}$.  
  It follows that 
  $dd^cu\wedge dd^cv\,\{0\}=C\cdot D$.  
  The tree measure of $u$ (resp. $v$) is supported on $\nu_C$
  (resp. $\nu_D$) and has mass $m(C)$ (resp. $m(D)$).  
  Whence 
  $\rho_u\cdot\rho_v=m(C)m(D)\,\nu_C\cdot\nu_D=C\cdot D$ 
  (see Section~\ref{sec-skewness}). This concludes the proof.
\end{proof}
\begin{proof}[Proof of Lemma~\ref{lem-log-sing}]
  Suppose $v=a\log\sum_i|f_i|^2$, $v'=a'\log\sum_i|f'_i|^2$,
  for holomorphic germs $f_i,f'_i\in R$, and $a,a'>0$. 
  Let $\pi$ be a resolution of singularities 
  of the curve $\{\prod_if_if'_i=0\}$, \ie 
  its total transform has normal crossings.
  It suffices to prove that $\pi^*(v-v')$ is locally bounded at any
  point on $\pi^{-1}\{0\}$.
  
  Pick a point $p\in\pi^{-1}(0)$.
  First suppose $p$ lies at the intersection of two irreducible components 
  $E$ and $F$ of $\pi^{-1}\{0\}$. 
  Choose coordinates such that $E\cup F$ is equal to
  $\{z=0\}\cup\{w=0\}$.  As the curve 
  $\{\prod_i\pi^*f_i=0\}$
  has normal crossings, it is equal to $E\cup F$ locally.
  This means that for any $i$ we can write 
  $\pi^*f_i=z^cw^d\xi$ for some $c,d\ge0$ and a unit $\xi\in R$.
  We infer that the function $\pi^*v$ differs from 
  $\div_E(\pi^*v)\log|z|+\div_F(\pi^*v)\log|w|$ 
  by a bounded function. 
  The same holds for $v'$.  
  As the tree transforms of $v$ and $v'$ coincide,
  the values of the divisorial valuations 
  $\nu_E$ and $\nu_F$ on $v$ and $v'$ are equal.
  Hence $\pi^*(v-v')$ is bounded at $p$.
  
  A similar argument applies when $p$ is a smooth point of $\pi^{-1}(0)$.
\end{proof}
%
%%%%%%%%%%%%%%%%%%%%%%%%%%%%%%%%%%%%%%%%%%%%%%%%%%%%%%%%%%%%%%%%%%%%%
% 
\subsection{Lower bound}
Our aim is to prove that one inequality in~\eqref{eq-inter} holds without
any restriction on the psh functions.
\begin{Prop}\label{p-lower}
  Suppose $u,v$ are psh functions such that the wedge product
  $dd^cu\wedge dd^c v$ is admissible. Then
  \begin{equation}\label{eq-lower}
    dd^cu\wedge dd^cv\,\{0\}
    \ge\rho_u\cdot\rho_v
    =\iint_{\cV\times\cV}\mu\cdot\nu\,d\rho_u(\mu)d\rho_v(\nu),
  \end{equation}
  where $\rho_u$, $\rho_v$ are the tree measures of 
  $u$ and $v$, respectively.
\end{Prop}
Note that in particular, when the right hand side in~\eqref{eq-lower}
is infinite, the wedge product $dd^c u\wedge dd^cv$ cannot be
admissible.

In view of~\eqref{e910} and Theorem~\ref{chartreetrans} we deduce the
following and significantly weaker classical result (see~\cite[Chapter~III,
Corollary~3.7.9]{dem0}).
\begin{Cor}\label{C902}
  If the wedge product $dd^cu\wedge dd^cv$ is admissible, then its mass
  at the origin is bounded from below by the product
  of the Lelong numbers of $u$ and $v$ at the origin.
\end{Cor}
\begin{proof}[Proof of Proposition~\ref{p-lower}]
  Let $u_n$, $v_n$ be the Demailly approximants of $u$ and $v$,
  respectively (see Section~\ref{sec-demailly}). 
  Using~\eqref{e911} and Proposition~\ref{p-inter} we infer
  \begin{equation*}
    dd^cu\wedge dd^cv\,\{0\}
    \ge\limsup_{n\to\infty}\, dd^cu_n\wedge dd^cv_n\,\{0\}
    =\limsup_{n\to\infty}\rho_n\cdot\sigma_n,
  \end{equation*}
  where $\rho_n$ and $\sigma_n$ are the measures represented 
  by $u_n$ and $v_n$, respectively. 
  Now Proposition~\ref{P901} gives
  $\limsup\rho_n\cdot\sigma_n\ge\rho\cdot\sigma$, completing
  the proof.
\end{proof} 
%
%%%%%%%%%%%%%%%%%%%%%%%%%%%%%%%%%%%%%%%%%%%%%%%%%%%%%%%%%%%%%%%%%%%%%
% 
\subsection{The \Holder continuous case}
\begin{Thm}\label{p-upper}
  Suppose $\varphi$ is a psh weight for which $e^{\varphi}$ is \Holder
  continuous in a neighborhood of the origin. Then for any psh function
  $u$ we have
  \begin{equation}\label{e-upper}
    dd^cu\wedge dd^c\varphi\,\{0\}
    =\rho_u\cdot\rho_\varphi
    =\iint_{\cV\times\cV}\mu\cdot\nu\,d\rho_u(\mu)d\rho_\varphi(\nu),
  \end{equation}
  where $\rho_u$ and $\rho_\varphi$ are the measures on
  $\cV$ represented by $u$ and $\varphi$, respectively.
\end{Thm}  
An equivalent formulation is
\begin{Cor}\label{C101}
  Under the assumptions of Theorem~\ref{p-upper}, 
  the generalized Lelong number defined by $\varphi$ is an average 
  of valuations:
  \begin{equation*}
    \nu_\varphi(u)=\int_\cV\mu(u)\,d\rho_\varphi(\mu),
  \end{equation*}
  for any psh function $u$.
\end{Cor}
\begin{Cor}\label{C102}
  Two psh weights whose exponentials are \Holder continuous define the
  same generalized Lelong numbers iff they have the same tree
  transform, or, equivalently, iff they have the same tree measure,
  or, yet equivalently, iff their pullbacks by any modification have 
  the same Lelong numbers at any point on the exceptional divisor.
\end{Cor}
\begin{Remark}
  If $\varphi$ is a \emph{homogeneous} psh weight in coordinates
  $(x,y)$, \ie 
  $\varphi(x,y)=\Phi(|x|,|y|)=c^{-1}\Phi(|x|^c,|y|^c)$,
  for all $c>0$, then its tree measure $\rho_\varphi$ 
  is supported on the set of monomial valuations in $(x,y)$.
  In this case, Corollary~\ref{C101} implies that $\nu_\varphi$ 
  is an average of monomial valuations.
  This was proved by Rashkovskii~\cite[Corollary~1]{rash}
  in any dimension.
\end{Remark}
\begin{proof}[Proof of Theorem~\ref{p-upper}]
  Let $\varphi_n$ be the Demailly approximating sequence of $\varphi$
  as in Section~\ref{sec-demailly}. As $e^\varphi$ is 
  \Holder continuous,~\eqref{e802} gives
  \begin{equation}\label{e803}
    \varphi-\frac{C}{n}
    \le\varphi_n
    \le\left(1-\frac{2}{nc}\right)\varphi+C,
  \end{equation}
  near the origin. Proposition~\ref{p-comparison} then implies
  $\nu_{\varphi}(u)\ge\nu_{\varphi_n}(u)\ge(1-2/nc)\,\nu_\varphi(u)$,
  hence $dd^cu\wedge dd^c\varphi_n\{0\}\to dd^cu\wedge dd^c\varphi\{0\}$.
  
  On the other hand~\eqref{e803} also gives
  $g_{\varphi_n}\ge g_\varphi\ge(1-2/nc)g_{\varphi_n}$, where 
  $g_\varphi$ and $g_{\varphi_n}$ are the tree transforms of $\varphi$ and
  $\varphi_n$, respectively.
  Thus
  \begin{multline*}
    dd^cu\wedge dd^c\varphi\{0\}
    =\lim_{n\to\infty}dd^cu\wedge dd^c\varphi_n\{0\}
    =\lim_{n\to\infty}\rho_u\cdot\rho_{\varphi_n}=\\
    =\lim_{n\to\infty}\int_\cV g_{\varphi_n}\,d\rho_u
    =\int_\cV g_{\varphi}\,d\rho_u
    =\rho_u\cdot\rho_{\varphi}.
  \end{multline*}
  Here the second equality follows from Proposition~\ref{p-inter}
  and the fourth from dominated convergence.
\end{proof}
%
%
%%%%%%%%%%%%%%%%%%%%%%%%%%%%%%%%%%%%%%%%%%%%%%%%%%%%%%%%%%%%%%%%%%%%%%%%%%%
%
%

\end{document}